\documentclass[11pt]{article}

\usepackage[a4paper,margin=1in]{geometry}
\usepackage{amsmath,amssymb,amsfonts,amsthm}
\usepackage{graphicx,mathrsfs,array,enumitem}
\usepackage[svgnames]{xcolor}
\usepackage{hyperref,tikz,float,authblk}
\usepackage{booktabs, multirow, array}
\usepackage{algpseudocode,algorithm,mathtools}

\newcolumntype{C}[1]{>{\centering\arraybackslash}m{#1}} 

\usepackage{cite}

\newcommand{\R}{\mathbb{R}}
\newcommand{\Rn}{\mathbb{R}^N}
\newcommand{\tpose}{^\mathsf{T}}
\DeclareMathOperator*{\argmin}{argmin}
\DeclareMathOperator*{\argmax}{argmax}
\DeclareMathOperator{\dive}{div}

\theoremstyle{definition}
\newtheorem{defn}{Definition}[section]

\theoremstyle{remark}
\newtheorem{rmk}{Remark}

\title{Agent-Based Optimal Control for Image Processing}
\author[1]{Alessio Oliviero\thanks{alessio.oliviero@uniroma1.it}}
\author[1]{Simone Cacace\thanks{simone.cacace@uniroma1.it}}
\author[1]{Giuseppe Visconti\thanks{giuseppe.visconti@uniroma1.it}}
\affil[1]{Department of Mathematics, Sapienza University of Rome, Rome, Italy}

\date{\today}

\begin{document}
\maketitle

\begin{abstract}
We investigate the use of multi-agent systems to solve classical image processing tasks, such as colour quantization and segmentation. We frame the task as an optimal control problem, where the objective is to steer the multi-agent dynamics to obtain colour clusters that segment the image. To do so, we balance the total variation of the colour field and fidelity to the original image. The solution is obtained resorting to primal-dual splitting and the method of multipliers. Numerical experiments, implemented in parallel with CUDA, demonstrate the efficacy of the approach and its potential for high-dimensional data.
\end{abstract}


\section{Introduction}

Image segmentation and colour quantization~\cite{SharmaAggarwal2010,Cheng2001} are fundamental tasks in image processing, computer vision, and pattern recognition. Their goal is to simplify the visual representation of an image while preserving its essential structure and perceptual content.

Classical approaches to these problems include clustering-based methods such as
the k-means~\cite{Shan2018}, hierarchical clustering~\cite{Arbelaez2011}, and evolutionary or meta-heuristic techniques such as genetic algorithms~\cite{Han2017}.
While these methods have proven effective in many practical applications, they typically require the number of clusters or segments to be fixed a priori, which may not be known in advance. Moreover, their performance can be sensitive to initialization, noise, and outliers, limiting their robustness. These approaches belong to the class of \emph{unsupervised learning} techniques~\cite{Yu2010,Mittal2021}.

On the other hand, \text{supervised learning} methods have gained increasing attention in recent years, particularly through the use of deep and shallow neural networks~\cite{Liu2021,Hesamian2019,Ronneberger2015,Chen2018}.
Well-trained networks can achieve high segmentation accuracy even for complex and
heterogeneous images.
However, the effectiveness of these methods strongly depends on the quality and quantity of the available training data.

In recent years, there has been growing interest in alternative frameworks based on multi-agent systems~\cite{Albi2017,Gong2023,MotschTadmor2014}, that describe the evolution of data distributions through collective or agent-based dynamics. Such approaches, inspired by models of self-organization in biology~\cite{CuckerSmale2007} and society~\cite{HegselmannKrause2002}, interpret clustering and denoising as emergent phenomena arising from local interactions among simple agents. This dynamical viewpoint offers a flexible and adaptive perspective on image processing, where segmentation and quantization can emerge naturally from the underlying interaction rules.

In agent-based models, each datum -- here, each pixel or colour intensity -- is interpreted as an agent whose state evolves dynamically through interactions with its neighbours. The collective behaviour of the system then gives rise to emergent structures that naturally correspond to clusters or coherent regions. The main advantage of these dynamical approaches lies in their ability to adapt to the geometry of the data. Interactions can be defined not only in terms of feature similarity but also in terms of spatial proximity within the image domain.
Particularly relevant examples in this context are the approaches based on the Kuramoto model for coupled oscillators~\cite{Liu2016}, on the Hegselmann-Krause consensus dynamics~\cite{Herty2020,Cabini2025}, and the mean-shift algorithm~\cite{ComaniciuMeer2002}.

While the dynamical perspective of multi-agent systems is advantageous for image processing, their connection to the well-established framework of variational methods can be further exploited to enhance their performance and provide a more rigorous mathematical foundation. Our goal is therefore to introduce a novel unsupervised technique that combines the strengths of both paradigms. On one hand, purely local methods, such as Yaroslavsky's neighborhood filters~\cite{yaroslavsky2012digital}, operate on a pixel's immediate surroundings and, consequently, have a restricted capacity to capture the global structures within an image. On the other hand, variational methods offer a more powerful and holistic framework by reformulating image processing tasks as energy minimisation problems. This paradigm was notably established by the Mumford--Shah model for image segmentation~\cite{mumford_shah}, which simultaneously seeks to find a piecewise smooth approximation of an image and the boundaries that define the smooth regions. Building upon this foundation, the Rudin--Osher--Fatemi model introduced total variation (TV) regularisation, a technique that effectively removes noise while preserving sharp edges~\cite{ROF_1992}. The variational framework has since become a dominant force in the field, with the development of non-local filtering~\cite{dabov2009bm3d, gilboa2009nonlocal, buades2011non} and total generalised variation~\cite{bredies2010total}. Furthermore, modern primal-dual formulations have provided efficient algorithms for solving the underlying optimisation problems, even in complex non-convex cases~\cite{CP_2011, lee2019finite, chambolle2024stochastic}.

In this paper, we consider an image as a discrete set of agents located on a two-dimensional spatial domain. Each pixel is represented by a position $(x_i,y_i) \in \Omega \subset \mathbb{R}^2$ and an associated intensity or colour value $c_i \in [0,1]$ (in the case of grey-scale images). The state of the system at a given time $t$ is therefore described by the collection $\{(x_i, y_i, c_i(t))\}_{i=1}^N$, where $N$ denotes the total number of pixels or agents.

The interaction between agents is designed to model both photometric similarity and spatial proximity. In particular, each agent tends to align its state $c_i$ with the states of nearby agents whose intensities are similar, while the strength of this influence decays with spatial distance.

The challenge, however, lies in determining how strongly the agents should interact at each time and at each spatial scale. Fixed interaction parameters may lead either to excessive smoothing or to poor denoising. 
To overcome this limitation, we introduce an optimal control approach, 
in which the interaction parameters evolve in time according to a control law designed to balance fidelity to the original image and regularization effects. The cost functionals we consider are inspired by the principles of variational image processing, while the dynamics of the system are governed by an alignment model applied to the image pixels. This approach uniquely integrates the global perspective of variational methods with the localised, adaptive nature of multi-agent systems.

In previous models~\cite{Herty2020,Cabini2025}, the agents represent the pixels of the image and are allowed to move within the spatial domain, while their associated colour or intensity values remain fixed. 
At the end of the evolution, the final colour assigned to each region is obtained by averaging the intensities of the pixels that have converged to the same spatial cluster. In other words, clustering occurs in the spatial domain, and the chromatic information is only used a posteriori to label the resulting groups.

In contrast, in the present work the positions of the agents are fixed, and the dynamics take place in the space of intensity values. Each pixel interacts with its neighbours according to both tonal and spatial similarity, but its colour value $c_i(t)$ evolves in time as a consequence of these interactions.

The entire numerical scheme is implemented in a parallel computing architecture using CUDA to handle the high dimensionality inherent in image data. This parallelization is naturally suited to our model, as the evolution of each pixel can be computed for each time step simultaneously across all pixels. Consequently, the synchronization of the device is necessary only once per time step, making the approach highly efficient. Our results demonstrate that this optimal control framework, particularly with the total variation cost, provides a flexible and effective method for colour quantization and image segmentation, successfully navigating the trade-off between clustering and feature preservation.

The rest of this paper is organized as follows. In Section \ref{sec:model}, we introduce the multi-agent system, detailing the differential equations that govern the evolution of the image and defining the space of admissible control functions. Section \ref{sec:OCP} is dedicated to the first optimal control formulation. We define a cost functional based on the $L^2$ norm of the intensity gradient, establish the well-posedness of the problem, and derive the corresponding optimality system. We then present the direct-adjoint looping (DAL) algorithm used for its numerical approximation and discuss the results of numerical tests, highlighting both the strengths and weaknesses of this approach. To address these limitations, in Section \ref{sec:TV} we introduce an improved model based on total variation minimisation. We reformulate the problem as a saddle-point problem for a regularised functional and detail the Alternating Direction Method of Multipliers (ADMM) used to solve it. This section concludes with a series of numerical experiments comparing the performance of the two proposed algorithms, demonstrating the superiority of the total variation approach in preserving sharp edges. Finally, Section \ref{sec:conclusions} provides concluding remarks, summarising our contributions and suggesting possible future lines of research.

\section{Mathematical model}
\label{sec:model}

We consider an image as a discrete set of agents located on a two-dimensional spatial domain. Let $I$ be a grey–scale image composed of $N$ pixels. We represent it as a rectangular matrix with $N_R$ rows and $N_C$ columns, whose entries contain the pixels' colour intensity normalised in $[0,1]$. To keep consistency with Cartesian coordinates, the rows of the matrix are ordered from bottom to top. We associate to each pixel an agent
\[
p_i = [x_i, y_i, c_i]^\top, \qquad i=1,\ldots,N,
\]
where $(x_i, y_i)$ are the (discrete) spatial coordinates of the pixel and $c_i$ is its grey-level value. Hence, $p_1=[x_1,y_1,c_1]$ corresponds to the bottom-left pixel of $I$ and $p_N=[x_N,y_N,c_N]$ to the top-right one, see Figure \ref{fig: pixels}.
\begin{figure}[ht]
	\centering
	\begin{tikzpicture}[scale=0.7,>=stealth]
		
		\def\nc{6}   
		\def\nr{5}   
		
		\foreach \x/\y/\g in {
			0/0/20, 1/0/35, 2/0/50, 3/0/65, 4/0/80, 5/0/30,
			0/1/40, 1/1/55, 2/1/70, 3/1/25, 4/1/45, 5/1/60,
			0/2/75, 1/2/30, 2/2/15, 3/2/50, 4/2/65, 5/2/85,
			0/3/60, 1/3/80, 2/3/35, 3/3/20, 4/3/55, 5/3/70,
			0/4/25, 1/4/45, 2/4/60, 3/4/75, 4/4/35, 5/4/10
		}{
			\fill[black!\g] (\x,\y) rectangle ++(1,1);
		}
		
		\draw[step=1cm, black, thin] (0,0) grid (\nc,\nr);
		
		\draw[->, thick] (-0.4,0) -- (\nc+0.6,0) node[right] {$x$};
		\draw[->, thick] (0,-0.4) -- (0,\nr+0.6) node[above] {$y$};
		
		\node at (\nc/2,-0.7) {$N_C$ columns};
		\node[rotate=90] at (-0.8,\nr/2) {$N_R$ rows};
		
		\draw[red, ultra thick] (0,0) rectangle (1,1);
		\fill[red] (0.5,0.5) circle (1.3pt);
		\node[below left, red] at (0.0,0.0) {$p_1=[x_1,y_1,c_1]\tpose$};
		
		\draw[blue, ultra thick] (\nc-1,\nr-1) rectangle (\nc,\nr);
		\fill[blue] (\nc-0.5,\nr-0.5) circle (1.3pt);
		\node[above right, blue] at (\nc-0.5,\nr-0.05) {$p_N=[x_N,y_N,c_N]\tpose$};
		
		\draw[green!60!black, ultra thick] (2,2) rectangle (3,3);
		\fill[green!60!black] (2.5,2.5) circle (1.3pt);
		\draw[dashed, green!60!black] (2.5,2.5) -- (7,3.6);
		\node[align=left, green!60!black] at (9.7,3.8)
		{$p_i=[x_i,y_i,c_i]\tpose$\\
			$(x_i,y_i)$ spatial coordinates\\
			$c_i\in[0,1]$ gray intensity};
		
		\node[below] at (0.5,0) {$1$};
		\node[below] at (\nc-0.5,0) {$N_C$};
		\node[left] at (0,0.5) {$1$};
		\node[left] at (0,\nr-0.5) {$N_R$};
		
	\end{tikzpicture}
	\caption{Representation of a grayscale image as a discrete set of agents $p_i=[x_i,y_i,c_i]^T$ on a two-dimensional domain. Pixel $p_1$ corresponds to the lower-left corner, while $p_N$ corresponds to the upper-right corner.}
	\label{fig: pixels}
\end{figure}

The key idea is to keep the position of the agents fixed, thereby preserving the spatial correlation among them, and let the colour component evolve over a finite time horizon $[0,T]$ according to a system of coupled ordinary differential equations, which, for each $i=1,\ldots,N$, reads
\begin{equation}
	\label{eq:MAsystem}
	\begin{cases}
		\displaystyle \dot c_i(t)
		= \frac{1}{N} \sum_{j=1}^N
		\phi\bigl(r_{ij}^\varepsilon(t)\bigr) \bigl(c_j(t)-c_i(t)\bigr), \quad t \in (0,T], \\
		c_i(0) = c_i^0.
	\end{cases}
\end{equation}
Here $\phi:\mathbb{R}\to\mathbb{R}$ is a differentiable interaction kernel that weighs the influence of neighbouring agents, and the quantity
\begin{equation}
	\label{eq:ellipsoid}
	r_{ij}^\varepsilon(t)
	= \frac{\varepsilon_x(t)}{2}\bigl[(x_j-x_i)^2+(y_j-y_i)^2\bigr]
	+ \frac{\varepsilon_c(t)}{2}\bigl(c_j(t)-c_i(t)\bigr)^2
\end{equation}
measures a distance between $p_i$ and $p_j$, as sketched in Figure \ref{fig: ellipsoid}. 
\begin{figure}[ht]
	\centering
	\begin{tikzpicture}[scale=2]
		
		\draw[->, thick, gray] (0,0) -- (2.2,0) node[right, text=black]{$x$};
		\draw[->, thick, gray] (0,0) -- (-0.8,-0.8) node[below left, text=black]{$y$};
		\draw[->, thick, gray] (0,0) -- (0,1.8) node[above, text=black]{$c$};
		
		
		\draw[dashed, blue] (1.5,0) arc (0:180:1.5 and 0.4);
		
		\draw[dashed, blue] (0,1) arc (90:270:0.4 and 1);
		
		\fill[blue!15, opacity=0.7] (0,0) ellipse (1.5 and 1);
		
		\draw[blue, thick] (0,0) ellipse (1.5 and 1);
		
		\draw[blue, thick] (-1.5,0) arc (180:360:1.5 and 0.4);
		
		\draw[blue] (0,-1) arc (-90:90:0.4 and 1);
		
		\draw[dashed, red, thick] (0,0) -- (1.5,0) node[midway, below]{$\propto \varepsilon_x(t)^{-1/2}$};
		\draw[dashed, red, thick] (0,0) -- (-0.4,-0.4) node[midway, left]{$\propto \varepsilon_x(t)^{-1/2}$};
		\draw[dashed, red, thick] (0,0) -- (0,1) node[midway, left]{$\propto \varepsilon_c(t)^{-1/2}$};
		
		\draw[->, orange, thick] (0,0) -- (0.7, 0.4) node[midway, above]{$r_{ij}^\varepsilon(t)\ $};
		
		\filldraw[black] (0,0) circle (0.8pt) node[anchor=east]{$p_i$};
		
		\filldraw[black] (0.7, 0.4) circle (0.8pt) node[anchor=south west]{$p_j$};
		
	\end{tikzpicture}
	\caption{Schematic representation of the interaction domain induced by the anisotropic distance $r_{ij}^{\varepsilon}(t)$ in the augmented space $(x,y,c)$. The kernel $\varphi$ selects the agents $p_j$ lying inside a bounded ellipsoidal neighbourhood of $p_i$, with semi-axes proportional to $\bigl(\varepsilon_x(t)^{-1/2},\varepsilon_x(t)^{-1/2},\varepsilon_c(t)^{-1/2}\bigr)$.}
	\label{fig: ellipsoid}
\end{figure}

The time-dependent parameters $\varepsilon_x(\cdot)$ and $\varepsilon_c(\cdot)$ play the role of anisotropic scaling factors, respectively quantifying the relevance of spatial separation and colour disparity in the interaction. When $\varepsilon_x$ is large compared with $\varepsilon_c$, agents tend to interact mostly with spatial neighbours, leading to diffusion–like smoothing; conversely, a small $\varepsilon_x$ enhances long–range connections based on colour similarity.

For any fixed $t \ge 0$, the set of $p_j$ that effectively influence $p_i$ is determined by the support of $\phi$ applied to $r_{ij}^\varepsilon(t)$. When $\phi$ has compact support, the interaction domain is a bounded ellipsoid in $[0,1]^3$, centred at $p_i$ and with semi–axes proportional to $(\varepsilon_x(t)^{-1/2}, \varepsilon_x(t)^{-1/2}, \varepsilon_c(t)^{-1/2})$. The time–dependence of $(\varepsilon_x,\varepsilon_c)$ allows the size of this ellipsoid to adapt, letting agents that were initially disconnected become linked during the evolution, or vice versa. Their evolution can be either fixed, scheduled, or determined through an optimal control strategy, as presented in the following section. Clearly, not every choice of  $\varepsilon := ( \varepsilon_x,\varepsilon_c)$ is suitable for our model. For example, negative values would have no physical meaning in this context.

\begin{defn}
	\label{def:control_space}
	Let $E$ be a convex, compact subset of $[0,+\infty)\times[0,+\infty)$. We define the space of \textit{admissible controls}
	\begin{equation}
		\label{eq:control_space}
		\mathcal{E} := \left\{ \varepsilon : \R^+_0 \to E \quad \text{measurable} \right\}.
	\end{equation}
\end{defn}

Classical results for opinion dynamics models~\cite{MotschTadmor2014} ensure that, under mild assumptions on $\phi$, the solution converges to a finite number of clusters in finite time. These clusters correspond to distinct regions of $\Omega$ over which the grey intensity is constant: a colour quantization of the original image.

\section{The optimal control problem}
\label{sec:OCP}

In this section, we formulate a control problem with the goal of retrieving an optimal quantisation of the input image. The first step is to define a cost functional that penalises colour variation, while at the same time keeping memory of the original image.

Suppose for a moment that $c$ is the colour field of an image defined on the continuous domain $\Omega := [0,1]^2$ and that $c^0$ is the corresponding version of the input image. In this setting, a simple cost functional with the properties outlined above would be
\begin{equation}
	\label{eq:cost_continuous}
	\mathcal{J}(c) = \int_\Omega \frac12 |\nabla c(x)|^2 + \frac{\alpha}{2}\left(c(x)- c^0(x)\right)^2\,dx, \quad \alpha>0.
\end{equation}

Before defining the corresponding discrete cost, we adopt the following notation.

\begin{defn}
	\label{def : discrete_gradient}
	Let $c=[c_i]_{i=1}^N$ be the vectorised representation of an $N$-pixel image with $N_R$ rows and $N_C$ columns, as described in Section \ref{sec:model}. Let also $\Delta x := 1/N_C$ and $\Delta y := 1/N_R$. We define the discrete Jacobian matrix $\nabla^\# c \in \R^{N \times 2}$ as the matrix whose rows are
	$$
	\nabla^\# c_i = \left[
	\frac{c_{i+1}-c_{i-1}}{2\,\Delta x}, \displaystyle \frac{c_{i+N_C}-c_{i-N_C}}{2\,\Delta y}
	\right],
	$$
	with reflective boundary conditions
	\begin{equation}
		\label{eq: reflective_boundary}
		\begin{aligned}
			i+1\ &\leftarrow\ i-1 &&\text{if } i \equiv 0 \text{ mod } N_C, \\
			i-1\ &\leftarrow i+1 &&\text{if } i \equiv 1 \text{ mod } N_C, \\
			i+N_C\ &\leftarrow\ i-N_C &&\text{if } i > N - N_C, \\
			i-N_C\ &\leftarrow\ i+N_C &&\text{if } i \le N_C.
		\end{aligned}
	\end{equation}
	In particular, we denote
	$$
	\left| \nabla^\# c_i \right|^2 =  \left(\frac{c_{i+1}-c_{i-1}}{2\,\Delta x}\right)^2 + \left( \frac{c_{i+N_C}-c_{i-N_C}}{2\,\Delta y} \right)^2.
	$$ 
\end{defn}
Then, by analogy with the continuous case, in our discrete setting we formulate the following cost.
\begin{defn}
	\label{def:cost_functional}
	Let $c=[c_i]_{i=1}^N$ be the vectorised representation of an $N$-pixel image as in Definition \ref{def : discrete_gradient}, $c^0=[c_i^0]_{i=1}^N$, and $\alpha > 0$. We define the cost function $J: \Rn \to \R$ as
	\begin{equation}
		\label{eq:cost_functional}
		J(c) = \frac12 \Delta x \, \Delta y\, \sum_{i=1}^N \left[ \left|\nabla^\# c_i\right|^2 + \alpha \left(c_i-c_i^0\right)^2 \right].
	\end{equation}
\end{defn}
The first term in \eqref{eq:cost_functional} penalises fluctuations in the colour field, while the second one enforces fidelity to the reference picture $c^0$. When applied to the solution of \eqref{eq:MAsystem}-\eqref{eq:ellipsoid}, the previous definition leads to the following Mayer problem, which falls within the framework of classical optimal control theory of ordinary differential equations:
\begin{equation}
	\label{eq:OCP}
	\underset{\varepsilon \in \mathcal{E}}{\text{minimise}} \quad J(c(T; \varepsilon)) \quad \text{subject to \eqref{eq:MAsystem}-\eqref{eq:ellipsoid}}.
\end{equation}
The cost $J$ is implicitly dependent on $\varepsilon = ( \varepsilon_x,\varepsilon_c)$, as these parameters are what determines the evolution from the initial condition  $c^0$ to $c(T)$. In this light, the cost functional can also be seen as the cost associated with the \textit{control policy} $\varepsilon \in \mathcal{E}$ over the time interval $[0,t]$. With a slight abuse of notation, in Algorithm \ref{alg:dal} we will refer to this control-to-cost map as $J(\varepsilon)$. The weight $\alpha$ balances the two cost components and is assumed to be a given parameter of the problem. We remark that the choice of this functional is motivated by the fact that the two terms in $J$ produce opposite behaviours: for $\alpha=0$ the minimiser is just a constant (single colour) image, whereas for $\alpha\to+\infty$ the minimizer will result in the original image $I$. Since the multi-agent dynamics is ultimately driven by the control policy $\varepsilon$, affecting spatial and colour interactions, these extremal configurations correspond respectively to full (one to all) and zero interactions, leading on the one hand to colour consensus, on the other hand to a clusterisation in which each pixel is itself a single cluster. Hence, by tuning $\alpha$, we expect to indirectly adjust the number of clusters in the final image. This may appear to be analogous to choosing a priori a fixed number of clusters in classical algorithms like $k$-means, but there is a fundamental difference between the two approaches. Selecting the number of segments is a \textit{prescriptive} constraint, while $\alpha$ is a \textit{descriptive} parameter. With classical methods, the user imposes a structure on the data and this is a strong assumption that may not be supported by the data itself. With our model, instead, the number of clusters is an emergent property. This is a data-driven result, not a predefined guess.

\begin{rmk}[On the choice of $T$ and $E$]
	\label{rmk:parameters}
	Although the optimal control problem can be stated for general $T \in (0,+\infty)$ and $E \subset\joinrel\subset [0,+\infty)\times[0,+\infty)$, from a practical point of view not every choice is equivalent. One should ensure that the compact set $E$ is large enough to allow both global interactions, i.e. $\phi(r_{ij}^\varepsilon)>0$ for any $i,j = 1,\ldots,N$, and no interactions at all, i.e. $\phi(r_{ij}^\varepsilon)=0$ for all $i,j = 1,\ldots,N$. Moreover, $T$ should be large enough to enable the formation of clusters -- or, in the language of opinion dynamics, to let the agents reach consensus. At the same time, an excessively far time horizon would just add a computational burden in the numerical simulations without yielding any benefit.
\end{rmk}

\subsection{Numerical approximation}
\label{sec:numerics_dal}

To solve \eqref{eq:OCP} we resort to Pontryagin's Maximum Principle (PMP)~\cite{PBGM_1964}. From this point onwards, we will write $F(c(t), \varepsilon(t))$ to refer to the right-hand side of \eqref{eq:MAsystem} at time $t$. When it does not create ambiguity, we will omit the dependency on $t$ to keep the notation as simple as possible.

We introduce the adjoint state function $\lambda : (0,T) \to \Rn$ and use it as a Lagrange multiplier to encode the dynamical constraints \eqref{eq:MAsystem}--\eqref{eq:ellipsoid} into an augmented cost functional, writing the Hamiltonian
$$
\mathcal{H}(c,\lambda,\varepsilon) = \lambda\tpose \cdot F(c, \varepsilon) = \sum_{i=1}^N \lambda_i\, F_i(c,\varepsilon).
$$
In this way, the PMP optimality system takes the form
\begin{subequations}
	\label{eq:opt_system}
	\begin{align}
		& \dot{c}(t) = F(c(t),\varepsilon(t)), \label{eq:opt_system-a} \\
		& c(0)=c^0, \notag \\
		& \dot{\lambda}(t) = - \partial_c F(c(t),\varepsilon(t))\tpose \cdot \lambda(t), \label{eq:opt_system-b} \\
		& \lambda(T) = - \nabla J(c(T)), \notag\\
		& \lambda(t)\tpose \cdot F(c(t), \varepsilon(t)) = \max_{\epsilon \in E} \ \lambda(t)\tpose \cdot F(c(t), \epsilon) \label{eq:opt_system-c}
	\end{align}
\end{subequations}
for $t \in (0,T)$. In particular, straightforward calculations yield:
\begin{align*}
	& \frac{\partial F_i}{\partial c_k} =
	\left\{
	\begin{aligned}
		&\frac1N \left[ \varphi'(r_{ik}^\varepsilon)\, \varepsilon_c (c_k-c_i)^2 + \varphi(r_{ik}^\varepsilon) \right] &\quad \text{if } k \neq i, \\
		&-\frac{1}{N} \sum_{j=1}^N \left[ \varphi'(r_{ij}^\varepsilon)\, \varepsilon_c ( c_j - c_i)^2 + \varphi(r_{ij}^\varepsilon)\right] &\quad \text{if } k = i,
	\end{aligned}
	\right. \\
	& \frac{\partial F_i}{\partial \varepsilon_x} = \frac{1}{2N} \sum_{j=1}^N \varphi'(r_{ij}^\varepsilon)\, \left( (x_j-x_i)^2 + (y_j - y_i)^2 \right)\, ( c_j - c_i), \\
	& \frac{\partial F_i}{\partial \varepsilon_c} = \frac{1}{2N} \sum_{j=1}^N \varphi'(r_{ij}^\varepsilon)\, ( c_j - c_i)^3, \\
	&\frac{\partial J}{\partial c_k} = \Delta x\, \Delta y \left[ -\Delta^\# c_k + \alpha (c_k-c^0_k) \right],
\end{align*}
where 
\begin{equation}
	\label{eq: discrete_laplacian}
	\Delta^\# c_k := \frac{c_{k+2} - 2 c_k + c_{k-2}}{(2\, \Delta x)^2} + \frac{c_{k+2N_C} - 2 c_k + c_{k-2N_C}}{(2\, \Delta y)^2}
\end{equation}
is a ``wide'' 5-point discrete Laplacian, and the same reflective boundary conditions \eqref{eq: reflective_boundary} are applied.

System \eqref{eq:opt_system} can be approximated with a simple direct-adjoint looping (DAL) algorithm, which is essentially a steepest descent along the functional $J$. For our purposes, we add a projection $\Pi_E$ onto the convex set $E$ to bound the control at each step. Moreover, we employ a two-way Armijo backtracking line search method \cite{A_1966} to select the best descent step-size at each iteration, resulting in Algorithm \ref{alg:dal}.
\begin{algorithm}[h]
	\caption{DAL algorithm for system \eqref{eq:opt_system} with two-way Armijo backtracking line search}
	\label{alg:dal}
	\begin{algorithmic}[1]
		\Require Initial guess for the control $\varepsilon^{(0)}$, $\sigma^{(0)}>0$, $b,c \in (0,1)$
		\State $j \gets 0$
		\Repeat
		\State Integrate \eqref{eq:opt_system-a} forward in time to obtain $c^{(j)}$
		\State Using $c^{(j)}$ for the terminal condition, integrate \eqref{eq:opt_system-b} backwards to get $\lambda^{(j)}$ 
		\State $D^{(j)} \gets \partial_\varepsilon [F(c^{(j)},\varepsilon^{(j)})\tpose \cdot \lambda^{(j)}]$
		\State $k \gets 0, \quad \tau^{(0)} \gets \sigma^{(j)}$
		\If{ $J(\Pi_E(\varepsilon^{(j)} - \sigma^{(j)}\, D^{(j)} )) > J(\varepsilon^{(j)}) - c \,\sigma^{(j)} |D^{(j)}|^2$ }
		\While{ $J(\Pi_E(\varepsilon^{(j)} - \tau^{(k)}\, D^{(j)} ) ) > J(\varepsilon^{(j)}) - c \,\tau^{(k)} |D^{(j)}|^2$ }
		\State $\tau^{(k+1)} \gets b\, \tau^{(k)}$
		\State $k \gets k+1$
		\EndWhile
		\Else
		\While{ $J(\Pi_E(\varepsilon^{(j)} - \tau^{(k)}\, D^{(j)} )) \le J(\varepsilon^{(j)}) - c \,\tau^{(k)} |D^{(j)}|^2$ }
		\State $\tau^{(k+1)} \gets \tau^{(k)}/b$
		\State $k \gets k+1$
		\EndWhile
		\EndIf
		\State $\sigma^{(j+1)} \gets \tau^{(k)}$
		\State $\varepsilon^{(j+1)} \gets \Pi_E ( \varepsilon^{(j)} + \sigma^{(j+1)}\ D^{(j)} )$
		\State $j \gets j+1$
		\Until{convergence}
	\end{algorithmic}
\end{algorithm}

Due to the non-linearity of the ODE system \eqref{eq:MAsystem}, the control-to-state map $\varepsilon \mapsto c(\,\cdot\, ,T; \varepsilon)$ is not convex. As a consequence, we can only expect convergence to a local minimiser.

\begin{rmk}[Convergence criteria for the DAL algorithm]
	\label{rmk:convergence}
	Since the control takes values in the convex, compact set $E$, the variational inequality 
	\begin{equation} \label{eq:variational_constraint}
		\partial_\varepsilon \mathcal{H}(c^*(t),\lambda^*(t),\varepsilon^*(t)) \cdot (\epsilon-\varepsilon^*(t)) \leq 0,
	\end{equation}
	for almost every $s \in (0,T)$ and all $\epsilon \in E$, which is equivalent to \eqref{eq:opt_system-c}, yields the stop condition
	$$
	\partial_\varepsilon \left[\sum_{i=1}^N \lambda_i^{(j)}\, F_i \left( c^{(j)},\varepsilon^{(j)} \right) \right] \cdot (\epsilon - \varepsilon^{(j)}) \leq \eta \quad \text{for a.e. } t \in (0,T) \text{ and all } \epsilon \in E,
	$$
	for some fixed tolerance $0< \eta \ll 1$. This relation is simplified in the case of box constraints, i.e.~when $E = [\varepsilon_x^{\min},\varepsilon_x^{\max}] \times [\varepsilon_c^{\min},\varepsilon_c^{\max}]$, as shown in~\cite{Troltzsch2024, Cacace_Oliviero}. Alternatively, a simpler criterion based on stationarity of the cost function can also be employed:
	$$
	\frac{|J(\varepsilon^{(j)}) - J(\varepsilon^{(j-1)})|}{J(\varepsilon^{(j-1)})} < \eta.
	$$
\end{rmk}

\subsection{Numerical tests}
\label{sec:test_dal}
In this section, we present numerical simulations obtained by applying Algorithm~\ref{alg:dal}. For all experiments, we set $T = 125$ and discretise the time interval $[0,T] = [0,125]$ with a uniform step size $\Delta t = 0.25$. The time integration is performed with an explicit Euler scheme.

The compactly supported interaction kernel is given by the $C^2$ \mbox{Wendland} function~\cite{Wendland1995}
\begin{equation}
	\label{eq:kernel}
	\varphi(r) :=
	\begin{cases}
		(4r-1)(1-r)^4, & 0 \leq r \leq 1, \\
		0, & \text{otherwise}.
	\end{cases}
\end{equation}
This kind of interaction kernels are standard in the smoothed particle hydrodynamics (SPH) literature, as they bring numerous advantages over gaussian or spline functions \cite{Dehnen2012}. They are less subject to numerical noise and reduce the formation of artefacts. Moreover, being piecewise polynomial, they are computationally cheap -- whereas exponential functions can impact performance -- and their derivatives are explicit, still piecewise polynomial, and numerically stable.

In accordance with Remark~\ref{rmk:parameters}, we set $E = [2,1100]^2$ and initialise the control as $\varepsilon_x^{(0)} = \varepsilon_c^{(0)} \equiv 57$ at all discrete times. For the two-way Armijo line search, we choose $\sigma^{(0)} = 1$, $b = \tfrac{1}{2}$, and $c = 10^{-4}$. Convergence is assessed via the stationarity of the cost functional -- computed with a rectangular quadrature rule -- with tolerance $\eta = 10^{-10}$ (see Remark~\ref{rmk:convergence}).

The routine is implemented in CUDA~C and executed on an NVIDIA H100~NVL GPU, hosted on the \textit{Lagrange} server of the Department of Mathematics, Sapienza University of Rome. In this way, the computations are parallelised over the pixels, yielding a significant performance improvement.

\begin{figure}[!h]
	\centering
	\begin{tabular}{cc}
		\includegraphics[width=.3\textwidth]{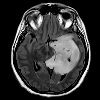} & \includegraphics[width=.4\textwidth]{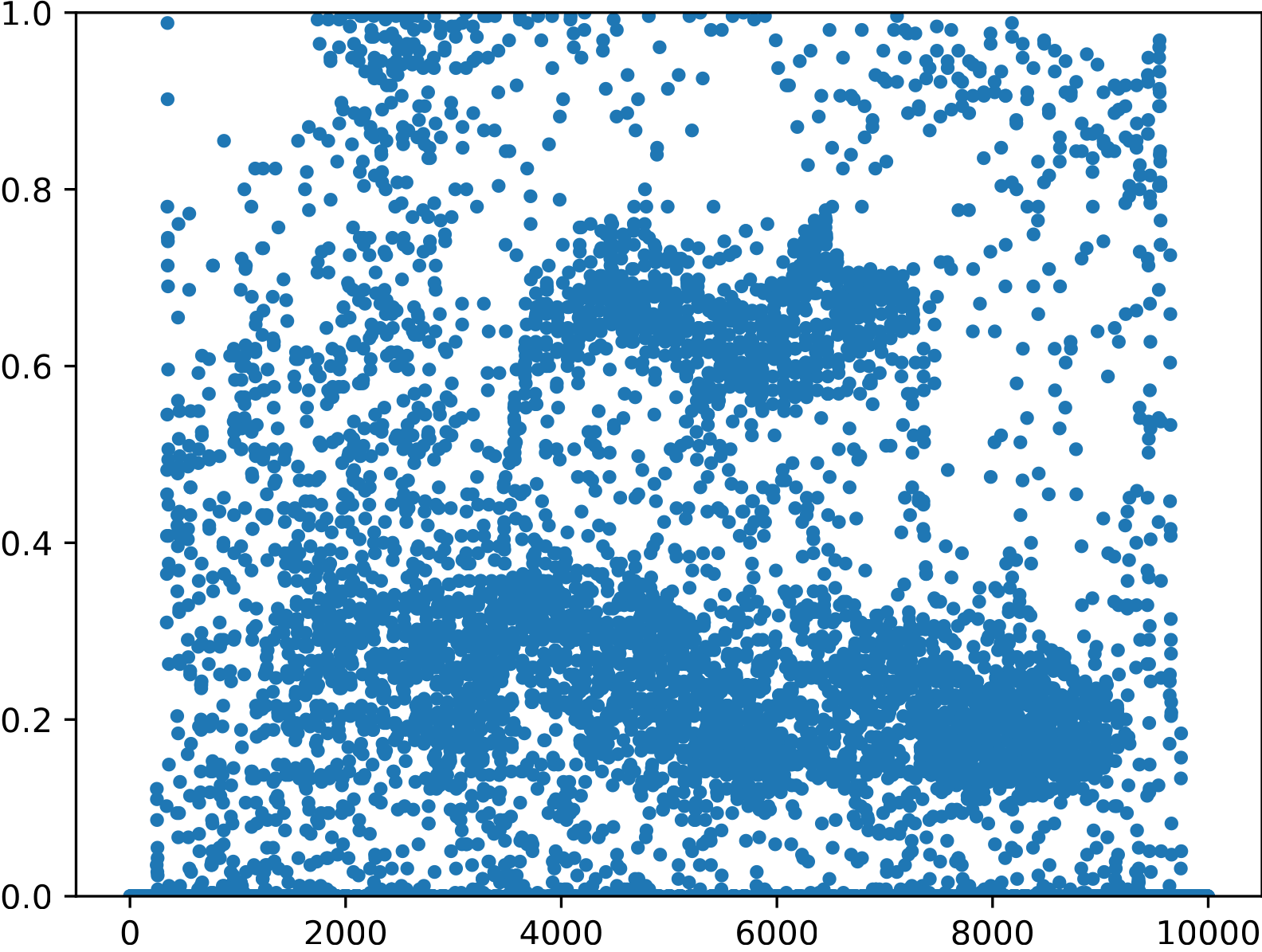}
	\end{tabular}
	\caption{Magnetic resonance image of a brain tumour, from \cite[Figure 1B]{mri_scan}. On the left side, the original picture; on the right, colour value $c_i \in [0,1]$ for every pixel $p_i,\ i=1,\ldots,N$.}
	\label{fig:mri}
\end{figure}

For the first experiment, we consider the image shown in the left panel of Figure~\ref{fig:mri}, which depicts a magnetic resonance image (MRI) of a brain tumour \cite[Figure~1B]{mri_scan}.

\begin{figure}[!h]
	\centering
	\begin{tabular}{lccc}
		\rotatebox{90}{$\qquad \alpha = 3000$} & 
		\includegraphics[width=.225\linewidth]{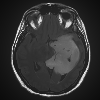} & 
		\includegraphics[width=.3\linewidth]{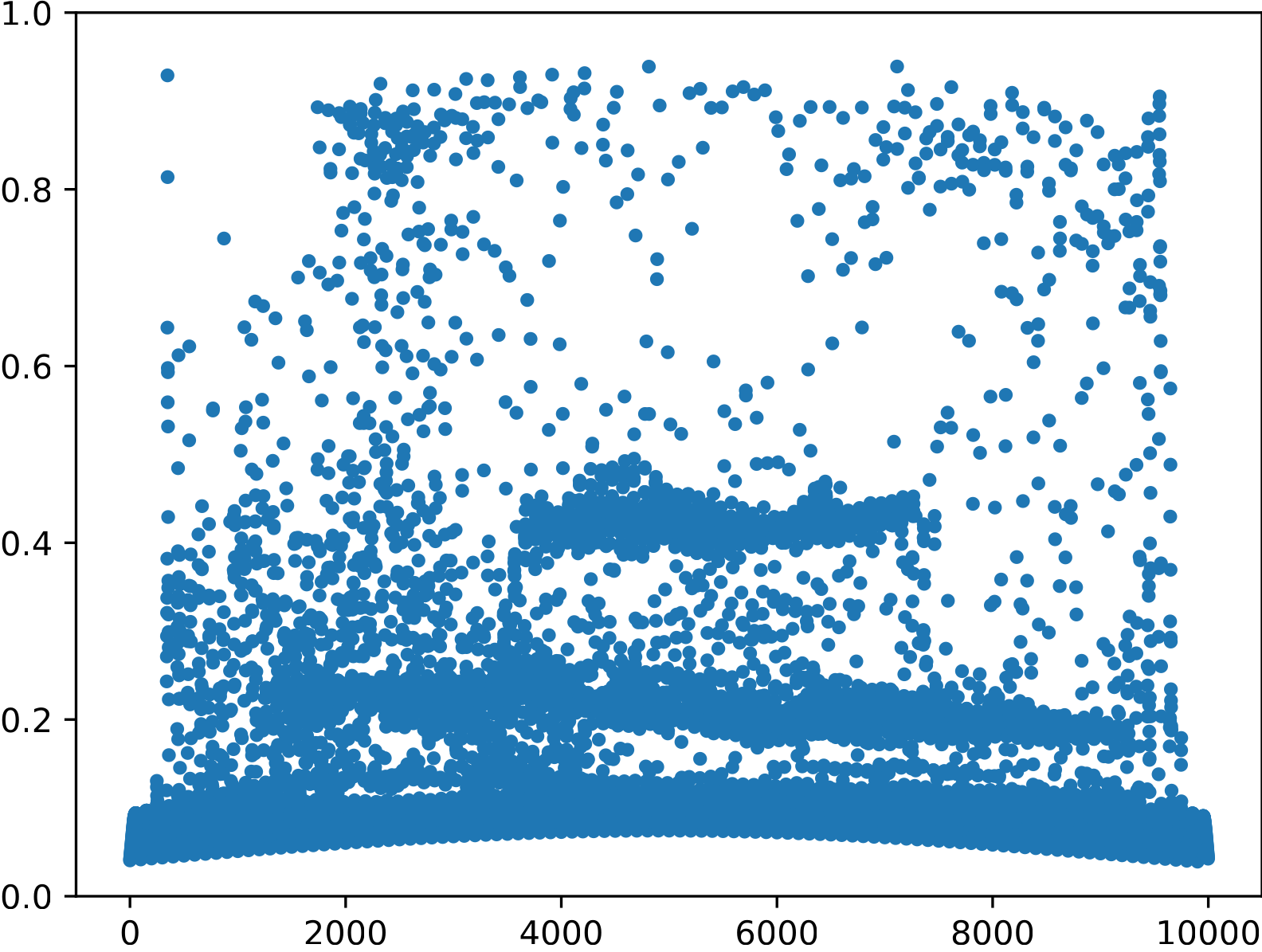} &
		\includegraphics[width=.3\linewidth]{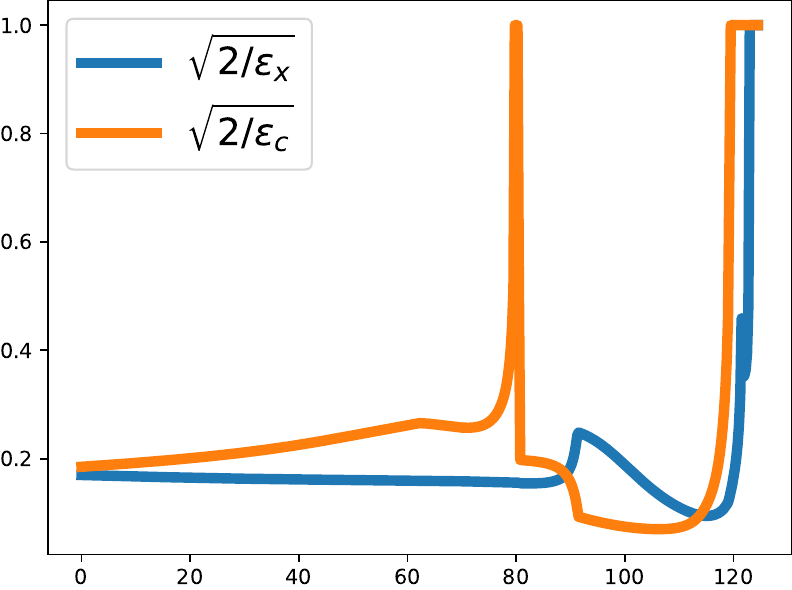} \\
		
		\rotatebox{90}{$\qquad \alpha = 1500$} & 
		\includegraphics[width=.225\linewidth]{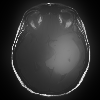} & 
		\includegraphics[width=.3\linewidth]{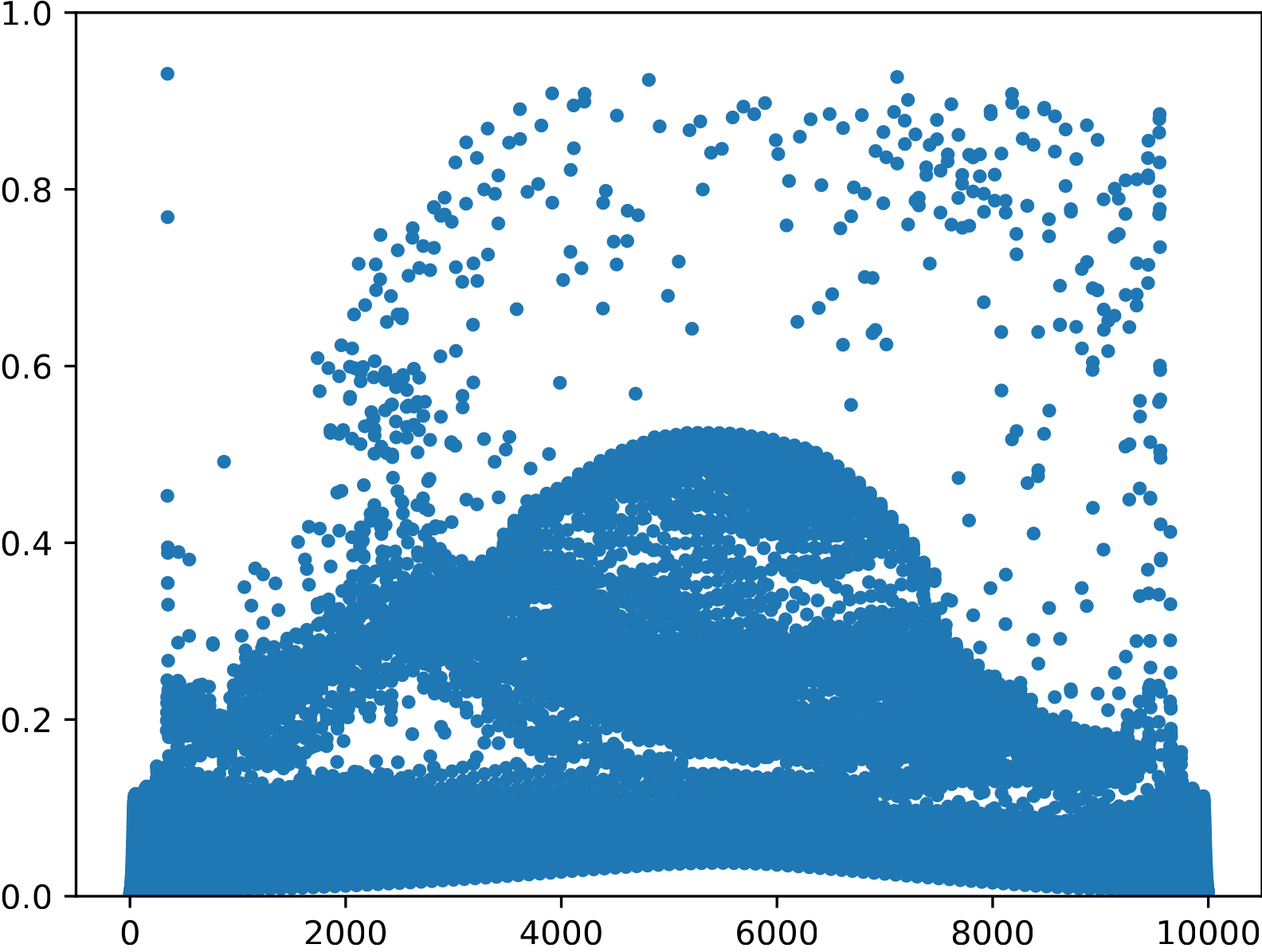} &
		\includegraphics[width=.3\linewidth]{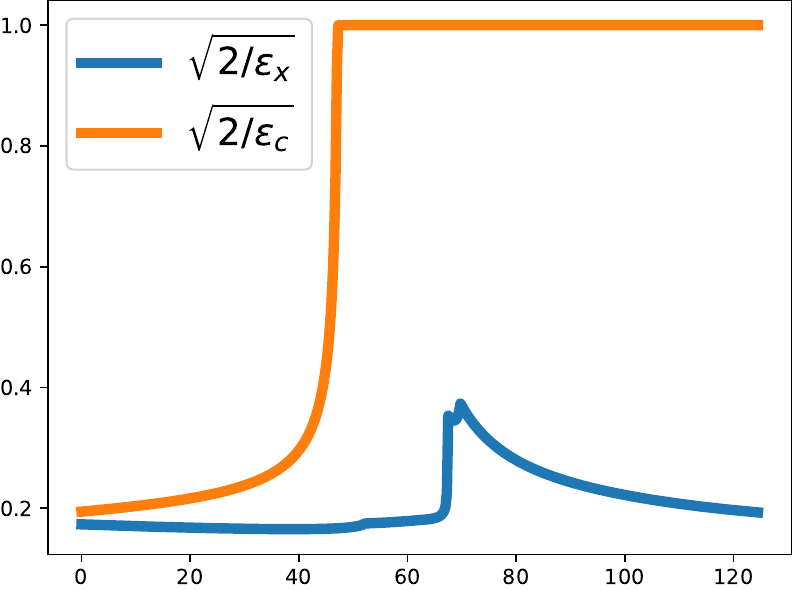} \\
		
		\rotatebox{90}{$\qquad \alpha = 600$} & 
		\includegraphics[width=.225\linewidth]{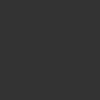} & 
		\includegraphics[width=.3\linewidth]{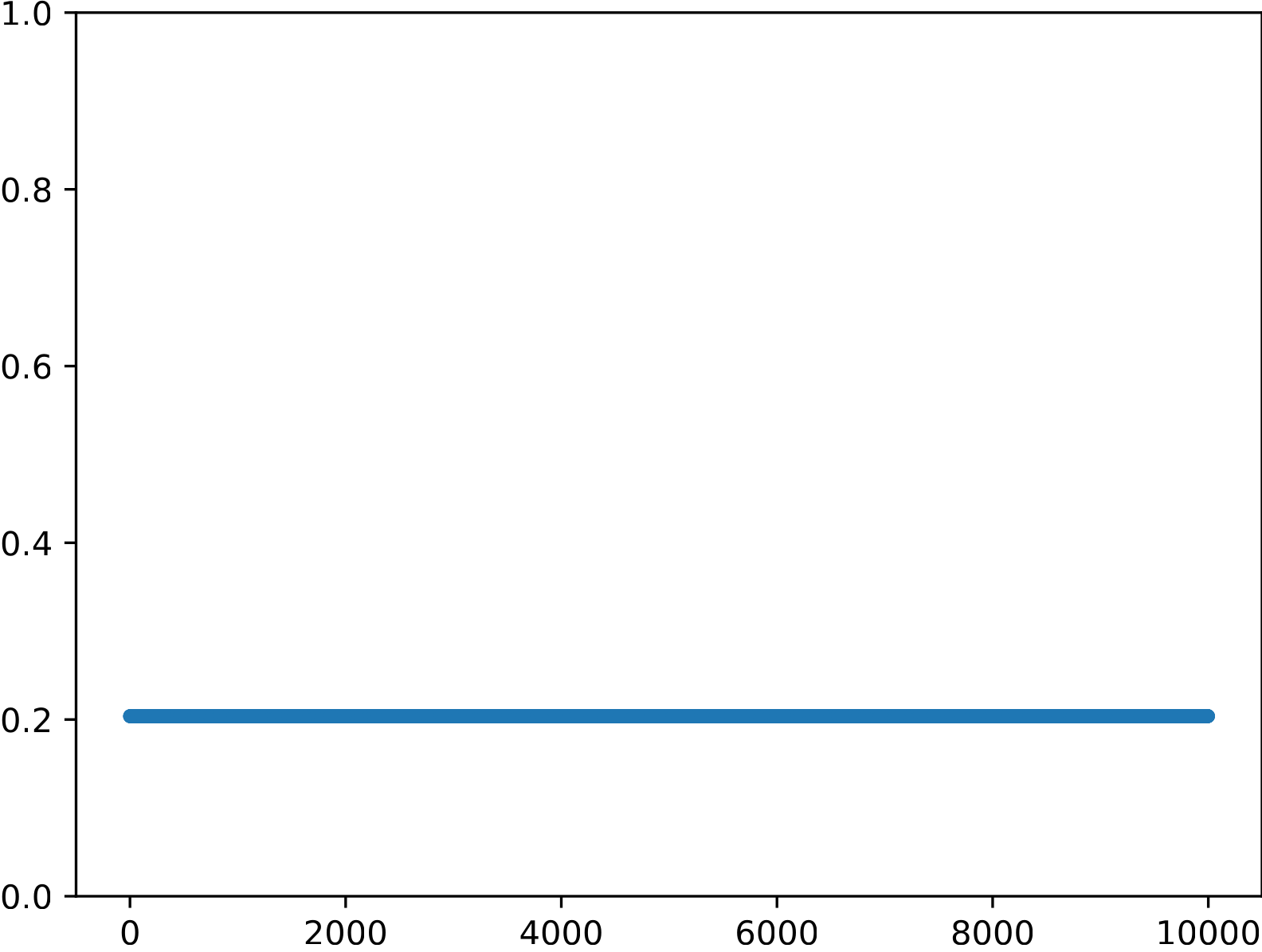} &
		\includegraphics[width=.3\linewidth]{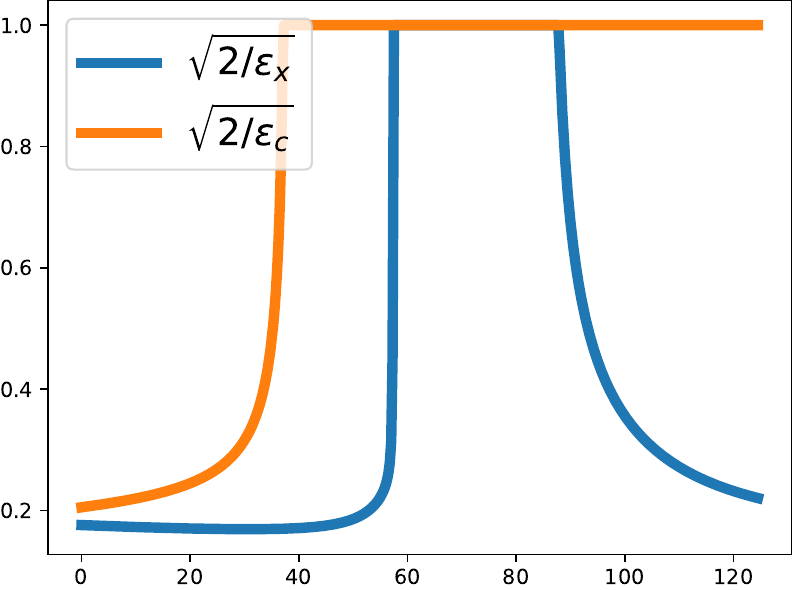}        
	\end{tabular}
	\caption{Output of Algorithm \ref{alg:dal} for different values of the fidelity parameter $\alpha$, using the image in Figure \ref{fig:mri} as input. On each row, from left to right: output image, grey level for every pixel, time evolution of the semi-axes of the ellipsoid \eqref{eq:ellipsoid}.}
	\label{fig:mri_dal}
\end{figure}

The colour distribution of the pixels is displayed in the right panel of the same figure. The output of Algorithm~\ref{alg:dal} for three different values of the fidelity parameter $\alpha > 0$ is reported in Figure~\ref{fig:mri_dal}. 
As expected, $\alpha$ has a decisive influence on the optimisation outcome. For large values of $\alpha$, the reconstruction closely follows the original image, and three distinct clusters emerge, corresponding to the background, healthy brain tissue, and tumour. For intermediate values, we observe a distinctive blurring effect caused by the penalisation of colour variation; consequently, edges are not well preserved. For smaller values of $\alpha$, the gradient term dominates, and global interactions drive the agents to consensus, producing a single cluster that represents the average colour of the image. This behaviour is reflected in the time evolution of the semi-axes of the interaction ellipsoid, which are plotted in the rightmost column of Figure~\ref{fig:mri_dal}. The smaller the value of $\alpha$, the more the control steers the system toward a regime of global interactions, as shown by the increase in $\varepsilon_x$ and $\varepsilon_c$ over time.

\begin{figure}[!h]
	\centering
	\begin{tabular}{cc}
		\includegraphics[width=.3\textwidth]{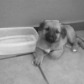} & \includegraphics[width=.4\textwidth]{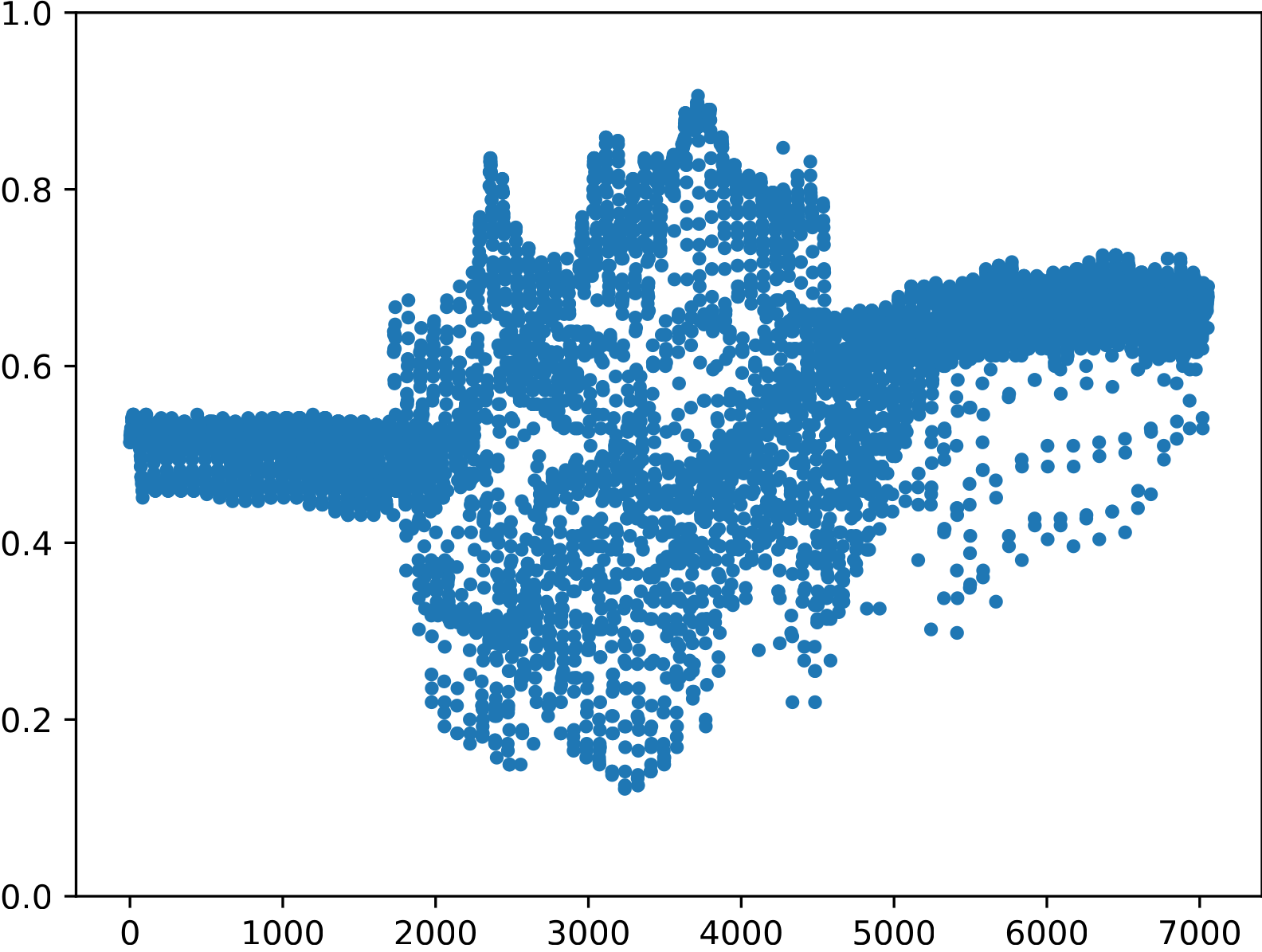}
	\end{tabular}
	\caption{Portrait of baby Namou, S. Cacace's dog. On the left side, the original picture; on the right, colour value $c_i \in [0,1]$ for every pixel $p_i,\ i=1,\ldots,N$.}
	\label{fig:namou}
\end{figure}

\begin{figure}[!h]
	\centering
	\begin{tabular}{lccc}
		\rotatebox{90}{$\qquad \alpha = 3000$} & 
		\includegraphics[width=.225\linewidth]{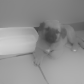} & 
		\includegraphics[width=.3\linewidth]{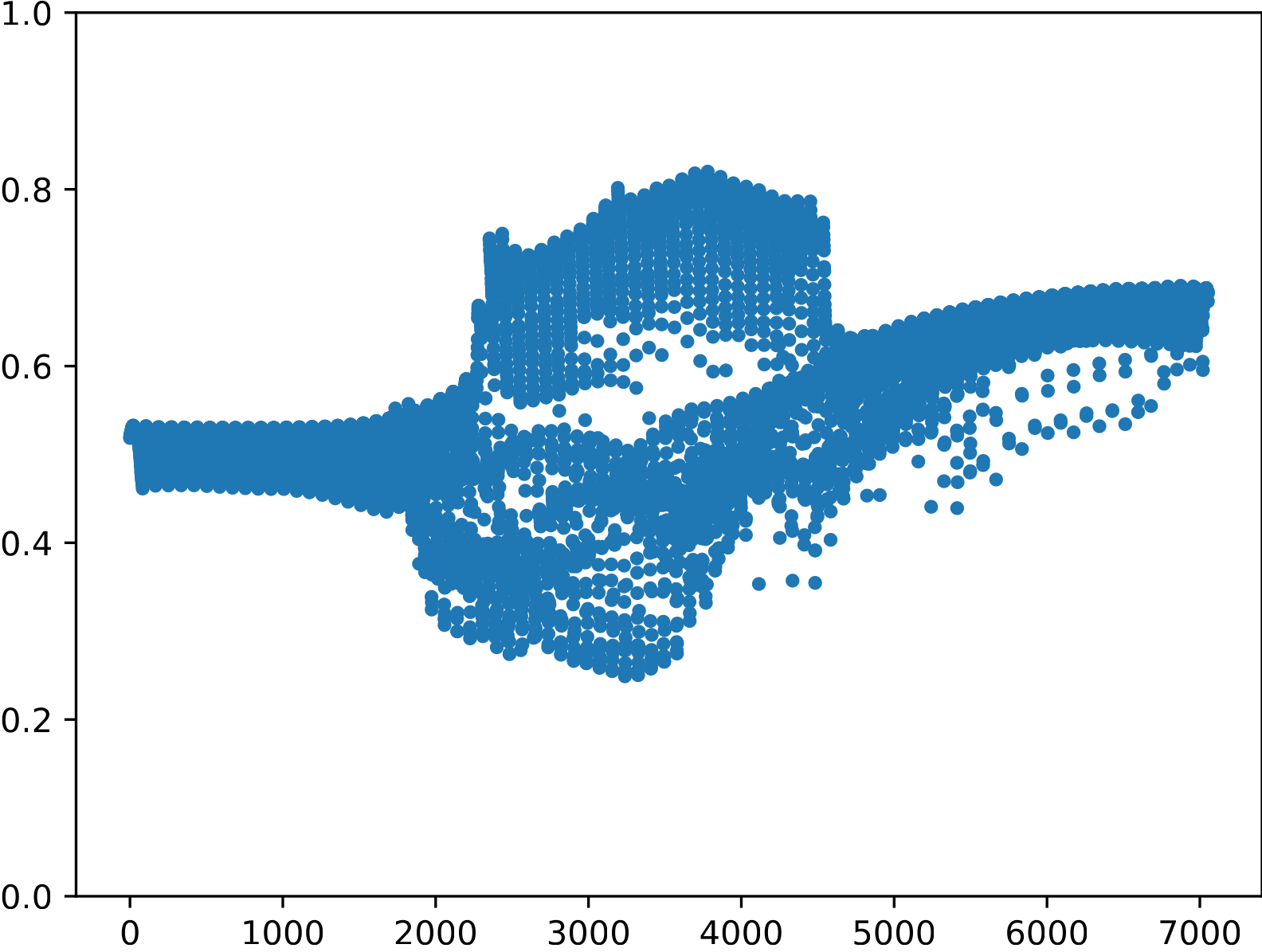} &
		\includegraphics[width=.3\linewidth]{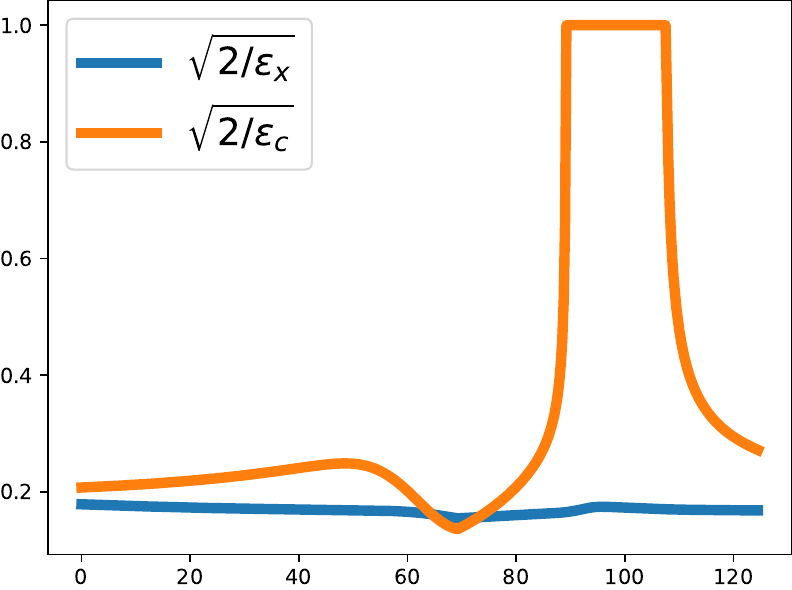} \\
		
		\rotatebox{90}{$\qquad \alpha = 1600$} & 
		\includegraphics[width=.225\linewidth]{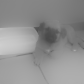} & 
		\includegraphics[width=.3\linewidth]{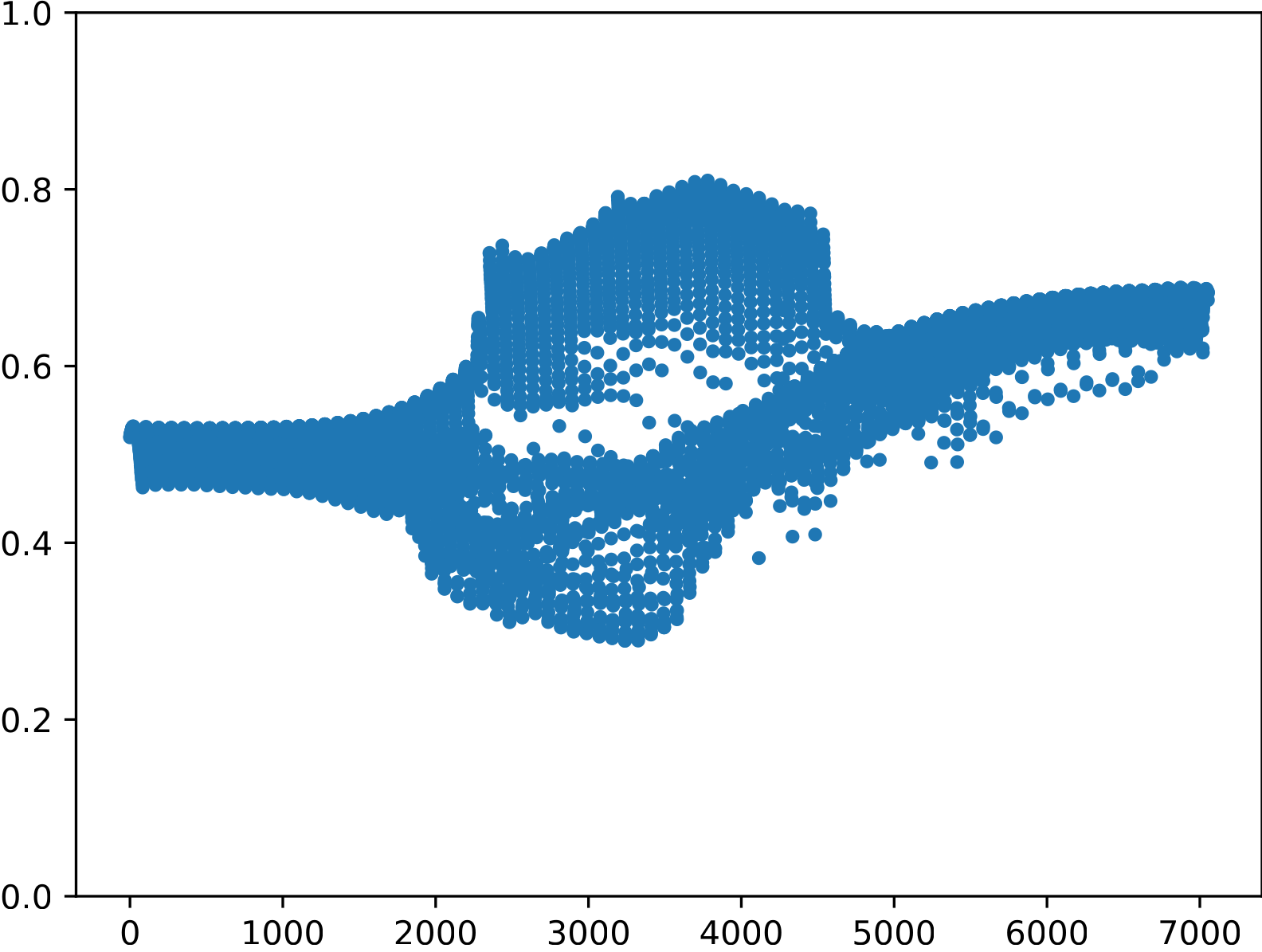} &
		\includegraphics[width=.3\linewidth]{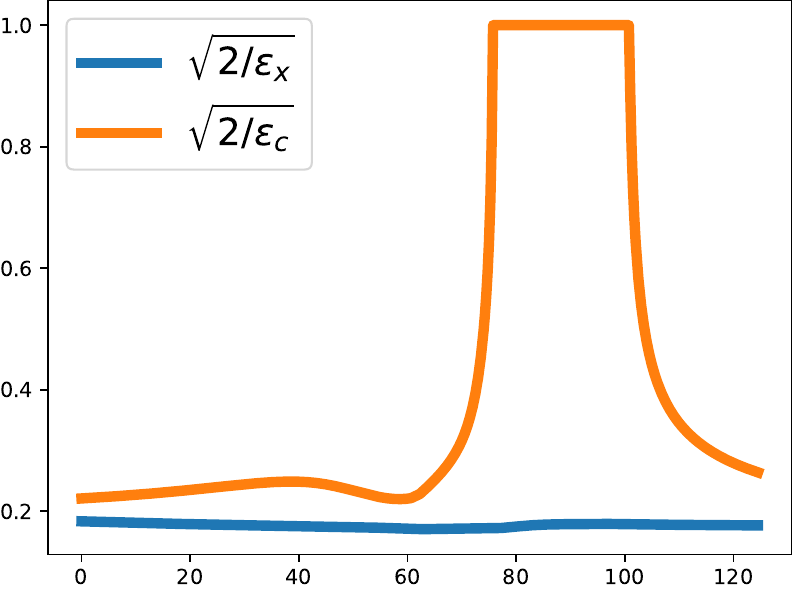} \\
		
		\rotatebox{90}{$\qquad \alpha = 600$} & 
		\includegraphics[width=.225\linewidth]{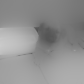} & 
		\includegraphics[width=.3\linewidth]{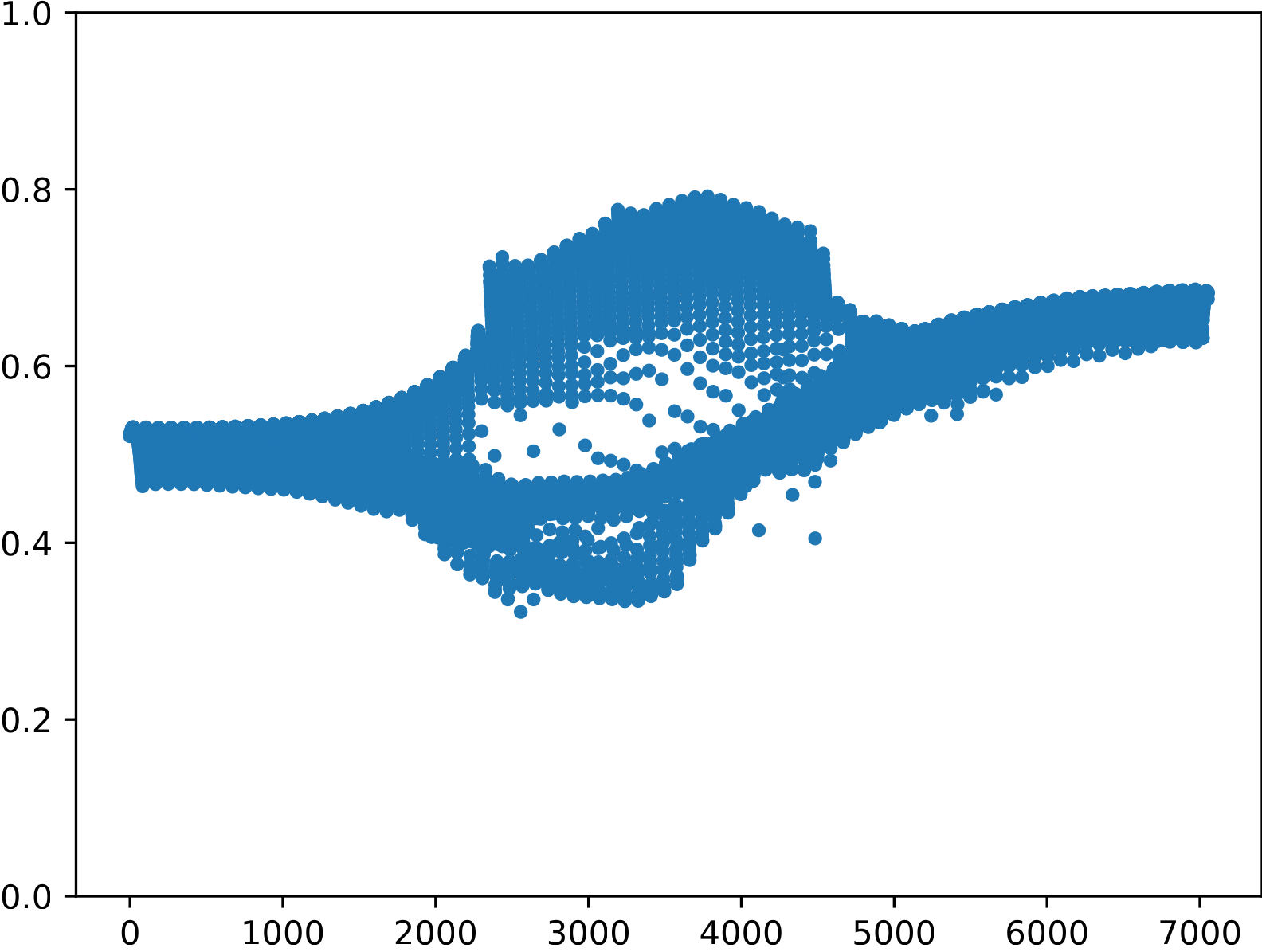} &
		\includegraphics[width=.3\linewidth]{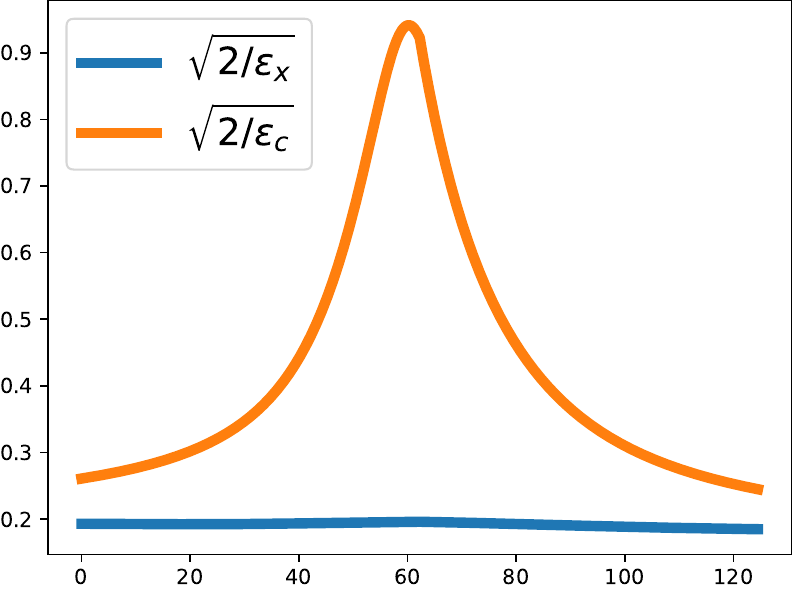}   
	\end{tabular}
	\caption{Output of Algorithm \ref{alg:dal} for different values of the fidelity parameter $\alpha$, using the image in Figure \ref{fig:namou} as input. On each row, from left to right: output image, grey level for every pixel, time evolution of the semi-axes of the ellipsoid \eqref{eq:ellipsoid}.}
	\label{fig:namou_dal}
\end{figure}

For a second experiment, we repeat the procedure on a more complex picture, shown in Figure~\ref{fig:namou}. 
The pixel distribution in the right panel reveals the absence of both purely black and purely white values. This makes the clustering task more challenging, as most agents are already concentrated within a narrow range of colour coordinates, despite representing different objects in the image. Figure~\ref{fig:namou_dal} illustrates the results: for large $\alpha$, the background wall and the floor are successfully flattened, but the bowl and the dog remain overly detailed. 
Reducing $\alpha$ improves the homogenisation of the bowl, but at the expense of blending the dog with the background. This time, the time evolution of the semi-axes of the ellipsoid \eqref{eq:ellipsoid} is less predictable. We can see from the rightmost column of Figure \ref{fig:namou_dal} that in all three cases the driving interaction is the one relative to colour, always reaching global (or close to global) interactions, whereas the space semi-axis is almost constant in time and always small. What changes is the moment in which global interactions are used (earlier for smaller values of $\alpha$) and the overall intensity over time, which is notably higher in the case $\alpha = 600$. 

We therefore conclude that Algorithm~\ref{alg:dal} is only effective for images with high contrast between structural elements, such as those arising in medical imaging, but less suited to natural images with smoother colour variations.

\section{Total variation minimisation}
\label{sec:TV}

Based on the previous tests, we introduce a different cost functional, with the aim of preserving sharp edges and mitigating the diffusion effect. To this end, we consider a fully discrete version of the functional from the Rudin--Osher--Fatemi (ROF) model~\cite{ROF_1992}
\begin{equation}
	\label{eq:ROF}
	\mathcal{K}(c) = \int_\Omega |\nabla c(x)| + \frac{\alpha}{2}(c(x)- c^0(x))^2\ dx,
\end{equation}
that is, in the notation of Definition \ref{def : discrete_gradient},
\begin{equation}
	\label{eq:cost_TV}
	K(c) = \Delta x \, \Delta y\, \sum_{i=1}^N \left[ \left|\nabla^\# c_i \right| + \frac{\alpha}{2}\left(c_i - c_i^0\right)^2 \right],
\end{equation}
with the same boundary conditions \eqref{eq: reflective_boundary}. In practice, we replace the quadratic term $\frac12 |\nabla^\# c_i|^2$ in \eqref{eq:cost_functional} with the total variation of $c$.
However, $K$ is not differentiable, and this prevents the immediate use of gradient-based numerical methods. To circumvent this difficulty, we resort to a primal--dual regularisation of the problem. 

\begin{defn}
	Let $v, \mu \in \R^{N\times 2}$ and $\alpha,\,\rho > 0$. With the notation of Definition \ref{def : discrete_gradient}, we define the regularised cost function
	\begin{multline}
		\label{eq:cost_tv_regular}
		L(c,v,\mu) = \Delta x \, \Delta y\, \sum_{i=1}^N \left[ |v_i| + \mu_i \cdot \left(v_i - \nabla^\# c_i \right) + \frac{\rho}{2} \left|v_i - \nabla^\# c_i \right|^2 \right. \\
		\left. + \frac{\alpha}{2}\left(c_i - c_i^0 \right)^2 \right].
	\end{multline}
	
\end{defn}
Considering the regularised cost and, like before, the implicit dependency of the solution of \eqref{eq:MAsystem}--\eqref{eq:ellipsoid} on the control $\varepsilon$, we obtain the optimal control problem
\begin{equation}
	\label{eq:OCP_TV}
	\underset{\substack{\varepsilon \in \mathcal{E}, \\ v \in \R^{N\times 2}}}{\text{minimise}} \ \underset{\mu \in \R^{N\times 2}}{\text{maximise}}  \quad L(c(T; \varepsilon), v, \mu) \quad \text{subject to \eqref{eq:MAsystem}--\eqref{eq:ellipsoid}}.
\end{equation}
In general, the saddle point we are looking for is not unique. In fact, when approaching a minimum point $(c^*,v^*)$ in the primal variables, the term $\sum_{i=1}^N \mu_i \cdot (v_i - \nabla^\# c_i)$ vanishes, so that any $\mu \in \R^{N\times 2}$ becomes feasible. However, this does not pose an issue for the numerical algorithms, which will converge to one of the possible saddle points for $L$.

\subsection{Numerical approximation}
To solve the saddle point problem \eqref{eq:OCP_TV}, we adopt the Alternating Direction Method of Multipliers (ADMM) \cite{Boyd2011}. This algorithm alternates a partial optimisation in each individual variable $c,v,\mu$, performed as outlined below, until a maximum number of iterations is reached, or until the gap between the subsequent decrease and increase in the cost function -- given, respectively, by the minimisation in the primal variables and the maximisation in the dual one -- is small enough. The routine is reported in Algorithm~\ref{alg:admm}. 

\begin{algorithm}[H]
	\caption{ADMM algorithm for the saddle point problem \eqref{eq:OCP_TV}}
	\label{alg:admm}
	\begin{algorithmic}[1]
		\Require Initial guesses for the control $\varepsilon^{(0)}$ and the dual variable $\mu^{(0)}$, $\gamma > 0$
		\State Integrate \eqref{eq:opt_system-a} forward in time to obtain $c^{(0)}$
		\State $v^{(0)} = \argmin_{v \in \R^{N\times 2}} L(c^{(0)}, v, \mu^{(0)})$
		\State $j \gets 1$
		\Repeat
		\State $\varepsilon^{(j)} = \argmin_{\varepsilon \in \mathcal{E}} L(c(T; \varepsilon), v^{(j-1)}, \mu^{(j-1)})$
		\State Integrate \eqref{eq:opt_system-a} forward in time to get $c^{(j)}$
		\State $v^{(j)} = \argmin_{v \in \R^{N\times 2}} L(c^{(j)}, v, \mu^{(j-1)})$
		\State $\mu^{(j)} = \argmax_{\mu \in \R^{N\times 2}} L(c^{(j)}, v^{(j)}, \mu)$
		\State $j \gets j+1$
		\Until{stop criterion is verified}
	\end{algorithmic}
\end{algorithm}

The key observation for ADMM is that, for any fixed $\mu$, the cost functional $L$ presents the same minimisers of
\begin{equation}
	\label{eq:complete_square}
	\begin{aligned}
		\tilde L (c,v) &:=  \Delta x \, \Delta y\, \sum_{i=1}^N \left[   \frac\rho2 \left|v_i - \nabla^\# c_i + \frac{\mu_i}{\rho} \right|^2 + |v_i| + \frac{\alpha}{2}\left(c_i - c_i^0 \right)^2 \right] \\
		& = L (c,v, \mu) + \frac{1}{2\rho} \Delta x \, \Delta y\, \sum_{i=1}^N |\mu_i|^2 .
	\end{aligned}
\end{equation}

\noindent As a consequence, for the minimisation in the $c$ variable we can adapt Algorithm~\ref{alg:dal}, substituting the final cost $J$ in the augmented Lagrangian formulation with the new $\tilde L$, assuming for the moment that $v$ is also fixed. With this procedure, we get to an optimality system that is the same as \eqref{eq:opt_system}, except for equation \eqref{eq:opt_system-b}. After some computations, we get
\begin{equation} \label{eq: admm_terminal}
	\frac{\partial \tilde L}{\partial c} = \Delta x\,\Delta y \left[ - \rho\, \Delta^\# c + \dive^\# (\mu + \rho\, v) + \alpha (c-c^0) \right],
\end{equation}
where the discrete divergence is defined as
\begin{equation}
	\label{eq: discrete_divergence}
	(\dive^\# w)_i := \frac{w_{i+1,1} - w_{i-1,1}}{2\, \Delta x} + \frac{w_{i+N_C,2} - w_{i-N_C,2}}{2\, \Delta y}
\end{equation}
for $w \in \R^{N\times 2}$, with the usual reflective boundary conditions \eqref{eq: reflective_boundary}. In particular, the terminal condition for the adjoint state $\lambda$ is now \eqref{eq: admm_terminal} computed at $c(T)$. It is worth noting that for $0 < \rho <\joinrel< 1$ the Laplacian term can become arbitrarily small. This means that now, by tuning the parameter $\rho$, we can mitigate diffusion, which is the cause of the blur effect in the tests presented in Section~\ref{sec:test_dal}.

For the minimisation in the $v$ variable, we observe that $v^*$ is a minimiser of $\tilde L (c, \, \cdot\,)$ if and only if $v^*_i$ minimises for all $i=1,\ldots,N$ the non-negative summand $ v_i \mapsto \frac{\rho}{2} \left|v_i - \nabla^\# c_i + \frac{\mu_i}{\rho} \right|^2 + |v_i|$. To keep the notation compact, consider the function $f(w) = \frac{\rho}{2}|w- a|^2 + |w|$, where $a$ is fixed. Then, if $|a| \le \frac{1}{\rho}$,
\begin{align*}
	f(w) &= \frac{\rho}{2} |w|^2 + \frac{\rho}{2} |a|^2 - \rho\, w \cdot a + |w| \\
	&\ge \frac{\rho}{2} |w|^2 + \frac{\rho}{2} |a|^2 - \rho |w| |a| + |w| \\
	&\ge \frac{\rho}{2} |a|^2 + ( 1- \rho|a|)\, |w| \ge \frac{\rho}{2} |a|^2 = f(0),
\end{align*}
for all $w$, implying $w^* = 0$. Conversely, if $|a| > \frac{1}{\rho}$, since $f$ is convex, for any $w \neq 0$ we can impose
\begin{equation}
	\label{eq:df}
	f'(w) = \rho (w-a) + \frac{w}{|w|} = 0.
\end{equation}
Moreover, in order to minimise $f$, it must be $w^* = C \frac{a}{|a|}$, for some $C > 0$. By plugging this expression into \eqref{eq:df}, we get
$$
\rho \left(C \frac{a}{|a|}-a\right) + \frac{a}{|a|} = 0,
$$
which yields $C = |a|-\frac{1}{\rho}$. Going back to the original problem, we obtain an explicit minimisation in the $v$ variable, the so-called shrinkage formula~\cite{Donoho1994,Boyd2011}:
\begin{equation}
	\label{eq:soft_threshold}
	\argmin_v L(c,v,\mu) = \argmin_v \tilde L(c,v) =
	\begin{cases}
		0 \quad &\text{if }\ |a| \le \frac{1}{\rho}, \\
		\displaystyle \left(|a|-\frac{1}{\rho} \right) \frac{a}{|a|} \quad &\text{otherwise,} 
	\end{cases}
\end{equation}
where $a = \nabla^\# c - \frac{1}{\rho}\mu$, and $| \cdot |$ is understood in the sense of Definition \ref{def : discrete_gradient}.

Finally, for the maximisation in the dual variable $\mu$, it is enough to perform a simple gradient ascent in the direction $v-\nabla^\# c$, since the cost functional $L$ is linear in $\mu$. Step 8 of Algorithm \ref{alg:admm} can be therefore implemented as
\begin{equation}
	\label{eq: dual_ascent}
	\mu^{(j)} \gets \mu^{(j-1)} + \gamma \left(v^{(j)} - \nabla^\# c^{(j)}\right).
\end{equation}

\begin{rmk}
	In the standard formulation of ADMM, the step-size $\gamma$ for the dual gradient ascent is set equal to the augmented Lagrangian penalty parameter $\rho$. This specific choice explicitly links the primal residual to the dual residual, so that the first-order optimality conditions are satisfied as the residuals converge to zero~\cite{Boyd2011}. However, in practice, it is often observed that tuning $\gamma$ independently from $\rho$ can lead to faster convergence~\cite{eckstein1992}. In our implementation, we set $\gamma = \rho$.
\end{rmk}

As already highlighted in Section~\ref{sec:numerics_dal}, we do not have convexity in the control-to-state map $\varepsilon \mapsto c(T; \varepsilon)$, so we can only expect convergence to a local minimum in the first variable. We also remark that, thanks to \eqref{eq:soft_threshold} and \eqref{eq: dual_ascent}, despite introducing two additional variables, the overall cost is that of just two nested loops. 

\subsection{Numerical tests}

In this section, we report a comparison between the outputs of Algorithm \ref{alg:dal} and Algorithm \ref{alg:admm} on the same test images used in Section \ref{sec:test_dal}. Additionally, we perform a third test on a more complex image. The discretisations and parameters used for the ADMM are exactly the same of Section \ref{sec:test_dal}, as well as the interaction kernel. The regularisation parameter is set to $\rho = 10^{-2}$ and the ascent step size is $\gamma = 10^{-2}$. Algorithm \ref{alg:admm} was also implemented in CUDA C and run on the same GPU as before.

\begin{figure}[!h]
	\centering
	\begin{tabular}{c c c}
		\textit{Original} & \textit{Algorithm \ref{alg:dal}} & \textit{Algorithm \ref{alg:admm}} \\
		\includegraphics[width=.3\textwidth]{MRI_mri.png} & 
		\includegraphics[width=.3\textwidth]{MRI_3000.png}  & 
		\includegraphics[width=.3\textwidth]{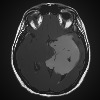} \\
		
		\includegraphics[width=.3\textwidth]{MRI_pixel_mri.png} &
		\includegraphics[width=.3\textwidth]{MRI_pixel_3000.png} &
		\includegraphics[width=.3\textwidth]{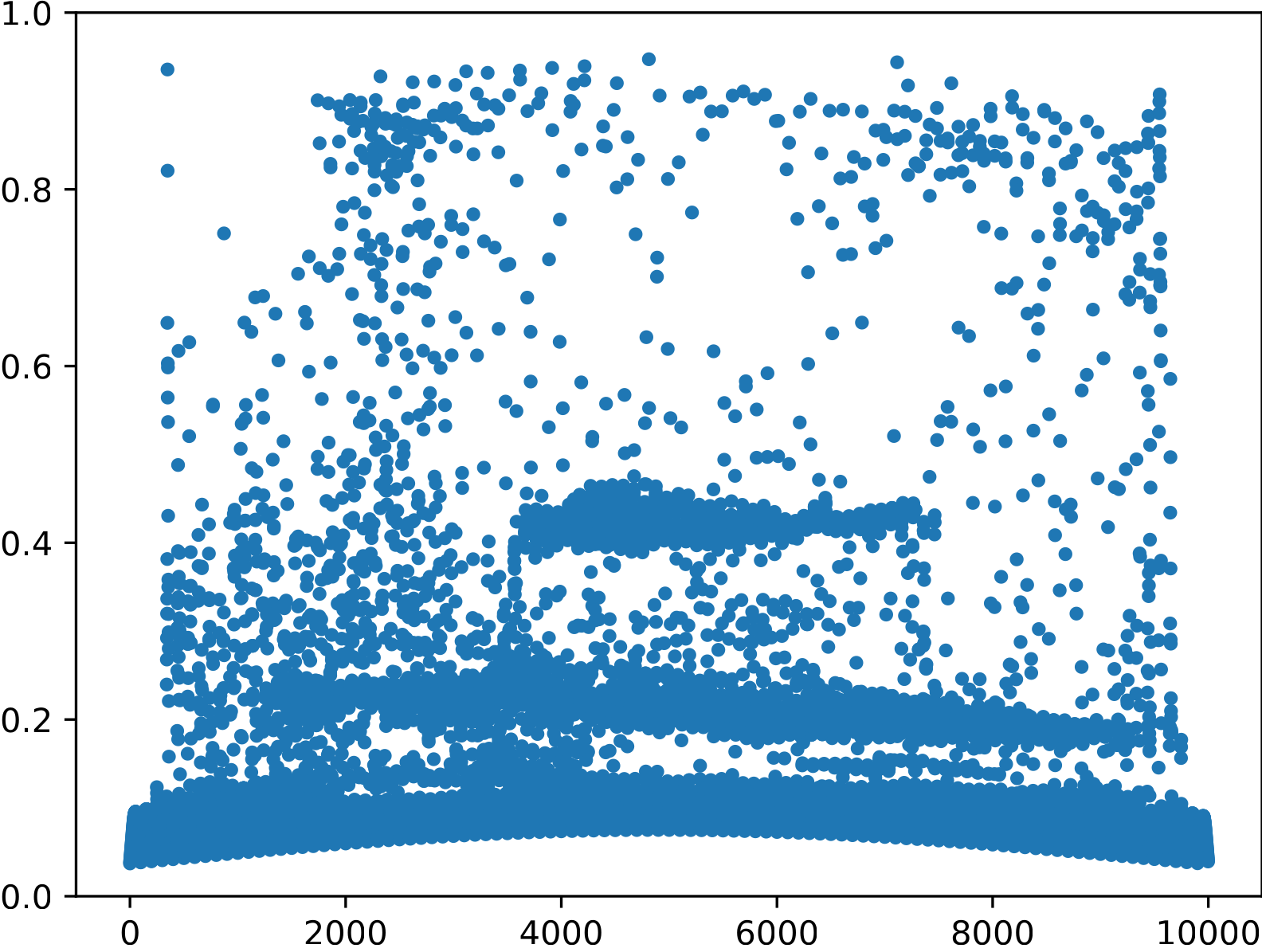} \\
		
		\ &
		\includegraphics[width=.3\textwidth]{MRI_E_3000.pdf} &
		\includegraphics[width=.3\textwidth]{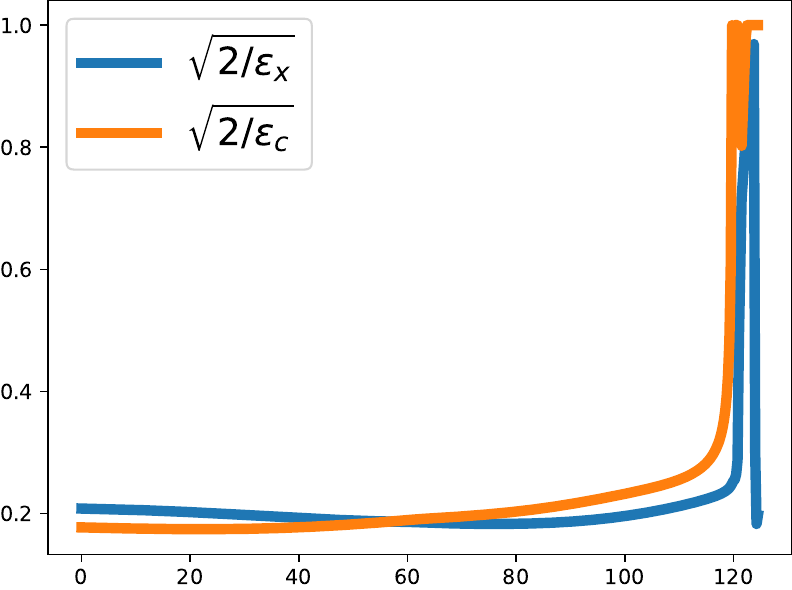}
	\end{tabular}
	\caption{Comparison of Algorithm \ref{alg:dal} and Algorithm \ref{alg:admm} on the test image in Figure \ref{fig:mri}. Output image on top, grey level for every pixel in the middle row, time evolution of the semi-axes of the ellipsoid \eqref{eq:ellipsoid} on the bottom.}
	\label{fig:mri_admm}
\end{figure}

The first test is carried out on the MRI scan in Figure \ref{fig:mri}, on which Algorithm \ref{alg:dal} works best. The results are shown in Figure \ref{fig:mri_admm}.
Despite finding qualitatively similar controls, reported on the bottom row, the healthy brain tissue and the tumour are depicted as more uniform clusters by Algorithm \ref{alg:admm}, without detailing them excessively, and at the same time without any blurring. In both cases, global interactions -- both in space and colour -- are used only at the end of the time interval.

\begin{figure}[!h]
	\centering
	\begin{tabular}{c c c}
		\textit{Original} & \textit{Algorithm \ref{alg:dal}} & \textit{Algorithm \ref{alg:admm}} \\
		\includegraphics[width=.3\textwidth]{NAM_namou.jpg} & 
		\includegraphics[width=.3\textwidth]{NAM_3000.png}  & 
		\includegraphics[width=.3\textwidth]{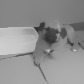} \\
		
		\includegraphics[width=.3\textwidth]{NAM_pixel_namou.png} &
		\includegraphics[width=.3\textwidth]{NAM_pixel_3000.png} &
		\includegraphics[width=.3\textwidth]{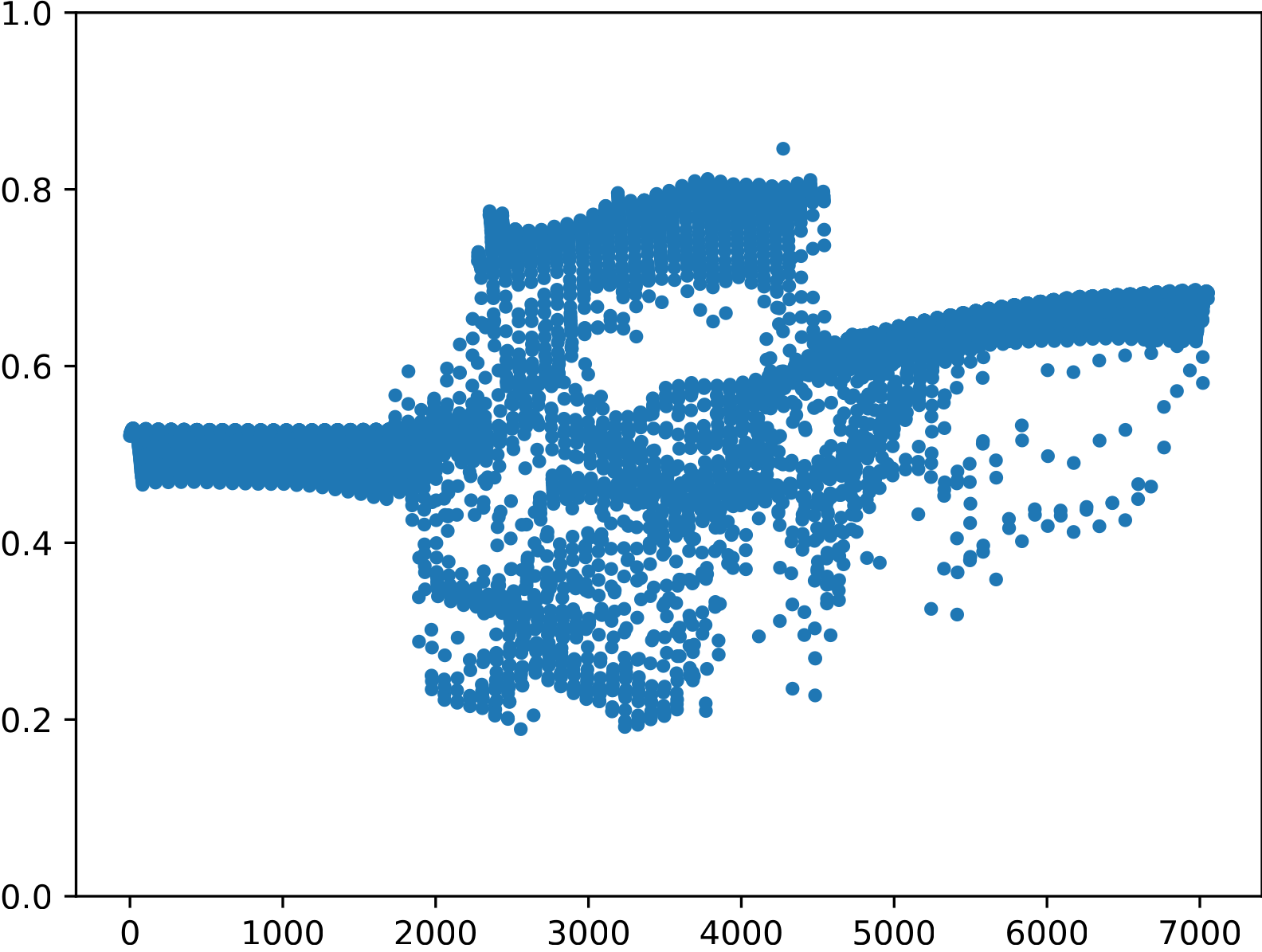} \\
		
		\ &
		\includegraphics[width=.3\textwidth]{NAM_E_3000.pdf} &
		\includegraphics[width=.3\textwidth]{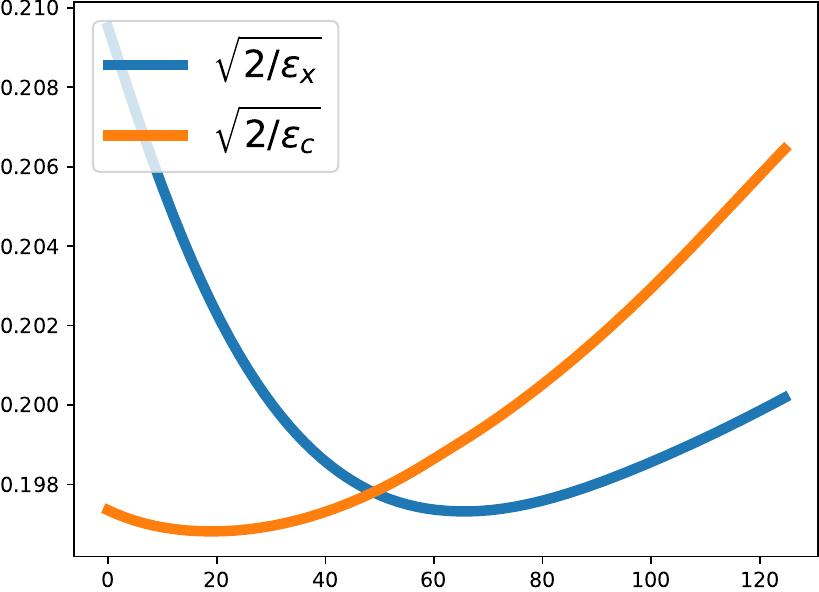}
	\end{tabular}
	\caption{Comparison of Algorithm \ref{alg:dal} and Algorithm \ref{alg:admm} on the test image in Figure \ref{fig:namou}. Output image on top, grey level for every pixel in the middle row, time evolution of the semi-axes of the ellipsoid \eqref{eq:ellipsoid} on the bottom.}
	\label{fig:namou_admm}
\end{figure}

Figure \ref{fig:namou_admm} illustrates the comparison using the dog picture in Figure \ref{fig:namou}. 

In this case, it is evident how Algorithm \ref{alg:admm} mitigates diffusion, unlike Algorithm \ref{alg:dal}, which tends to blur the image as already mentioned. Not only do the edges appear sharper in the top-right picture, but the segments themselves are more even. As a matter of fact, the semi-axes of the ellipsoid \eqref{eq:ellipsoid} are radically different. The ones found by Algorithm \ref{alg:admm} never involve global interactions (note the different scales in Figure \ref{fig:namou_admm}), which are clearly disadvantageous for this picture. This can be verified also by looking at the colour distributions in the middle row: the clusters created by Algorithm \ref{alg:admm} are more concentrated and fewer pixels exhibit an intermediate colour. 

\begin{figure}[!h]
	\centering
	\begin{tabular}{c c c}
		\textit{Original} & \textit{Algorithm \ref{alg:dal}} & \textit{Algorithm \ref{alg:admm}} \\
		\includegraphics[width=.3\textwidth]{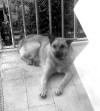} & 
		\includegraphics[width=.3\textwidth]{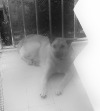}  & 
		\includegraphics[width=.3\textwidth]{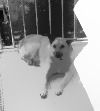} \\
		
		\includegraphics[width=.3\textwidth]{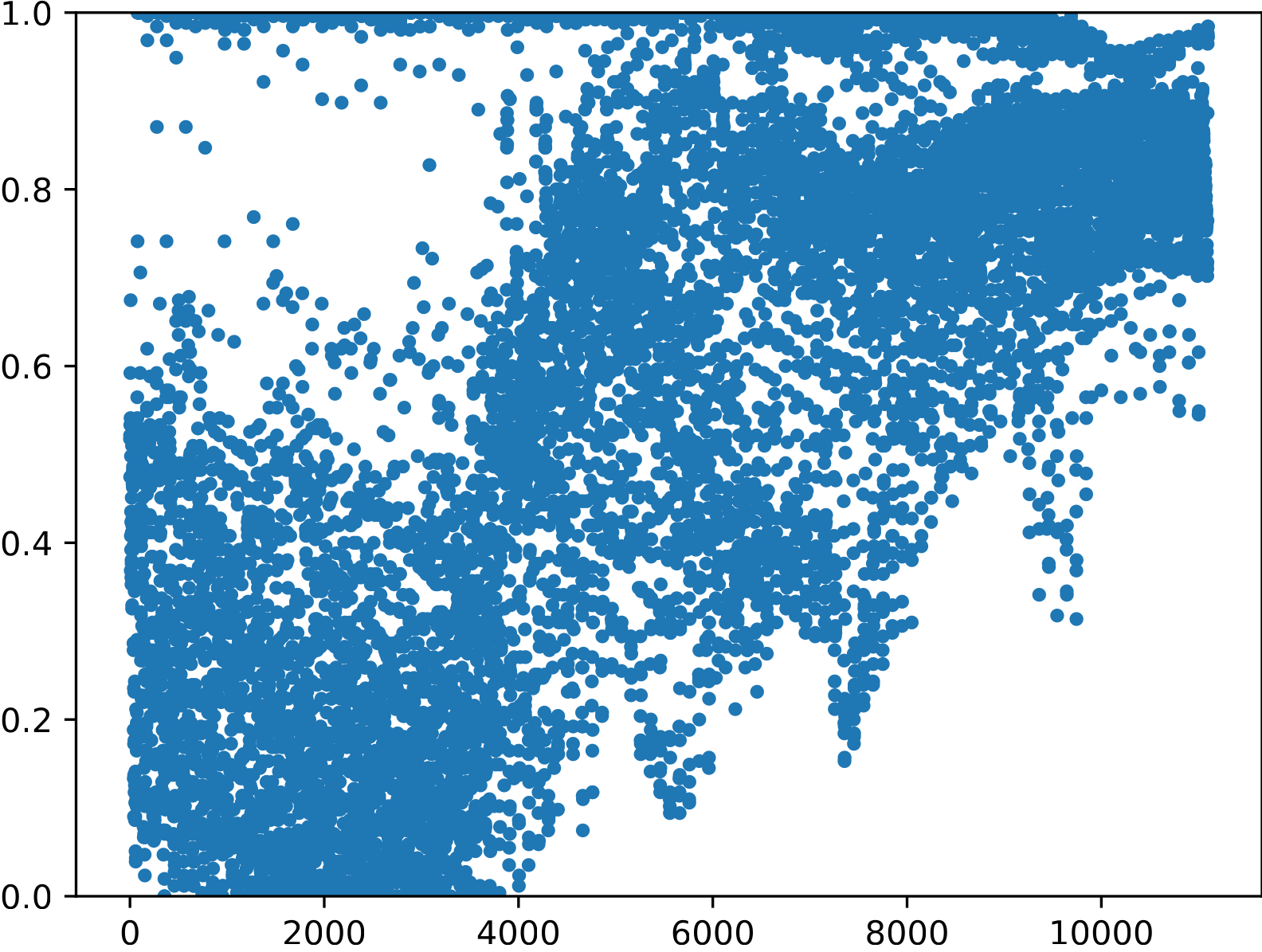} &
		\includegraphics[width=.3\textwidth]{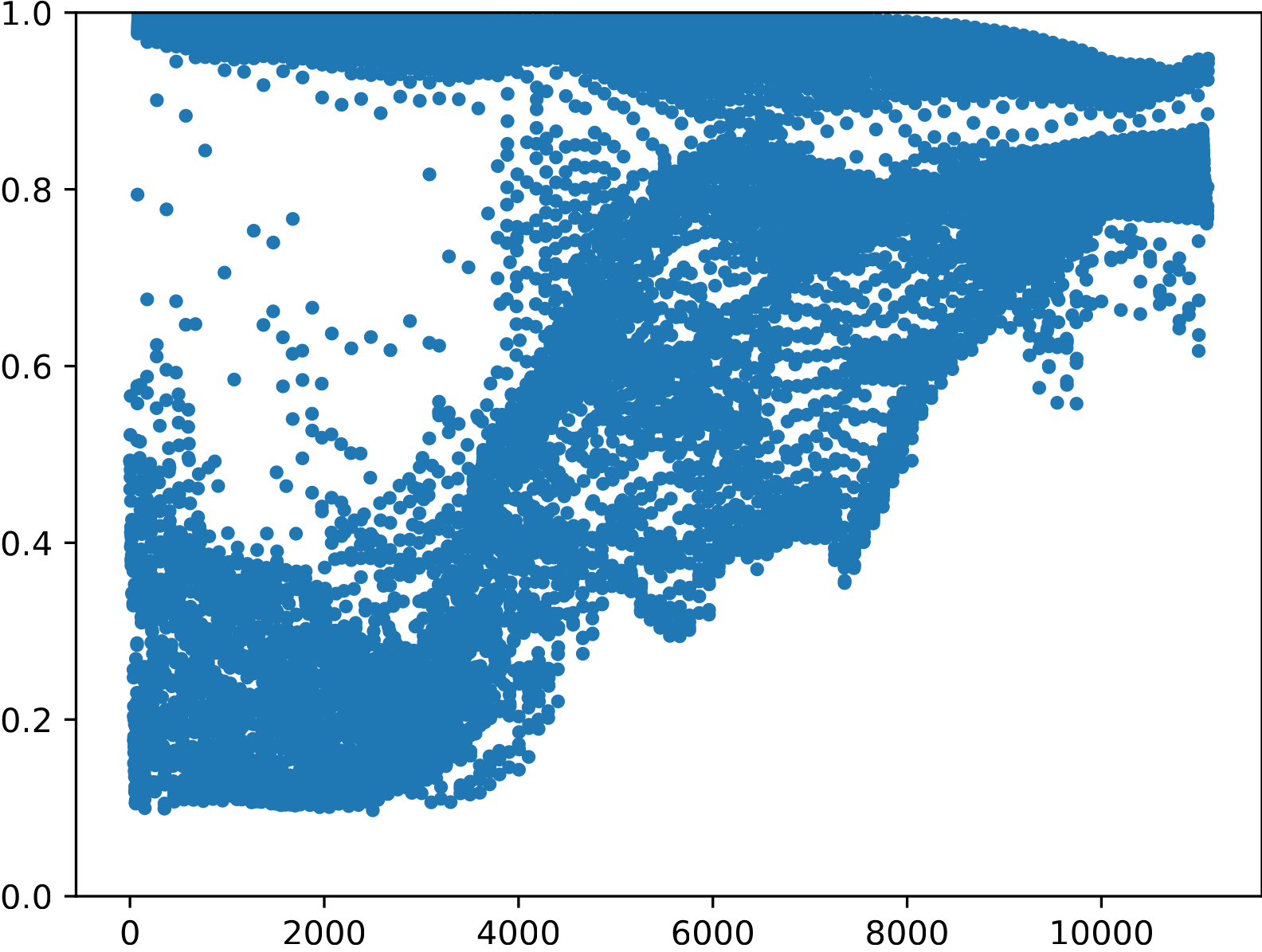} &
		\includegraphics[width=.3\textwidth]{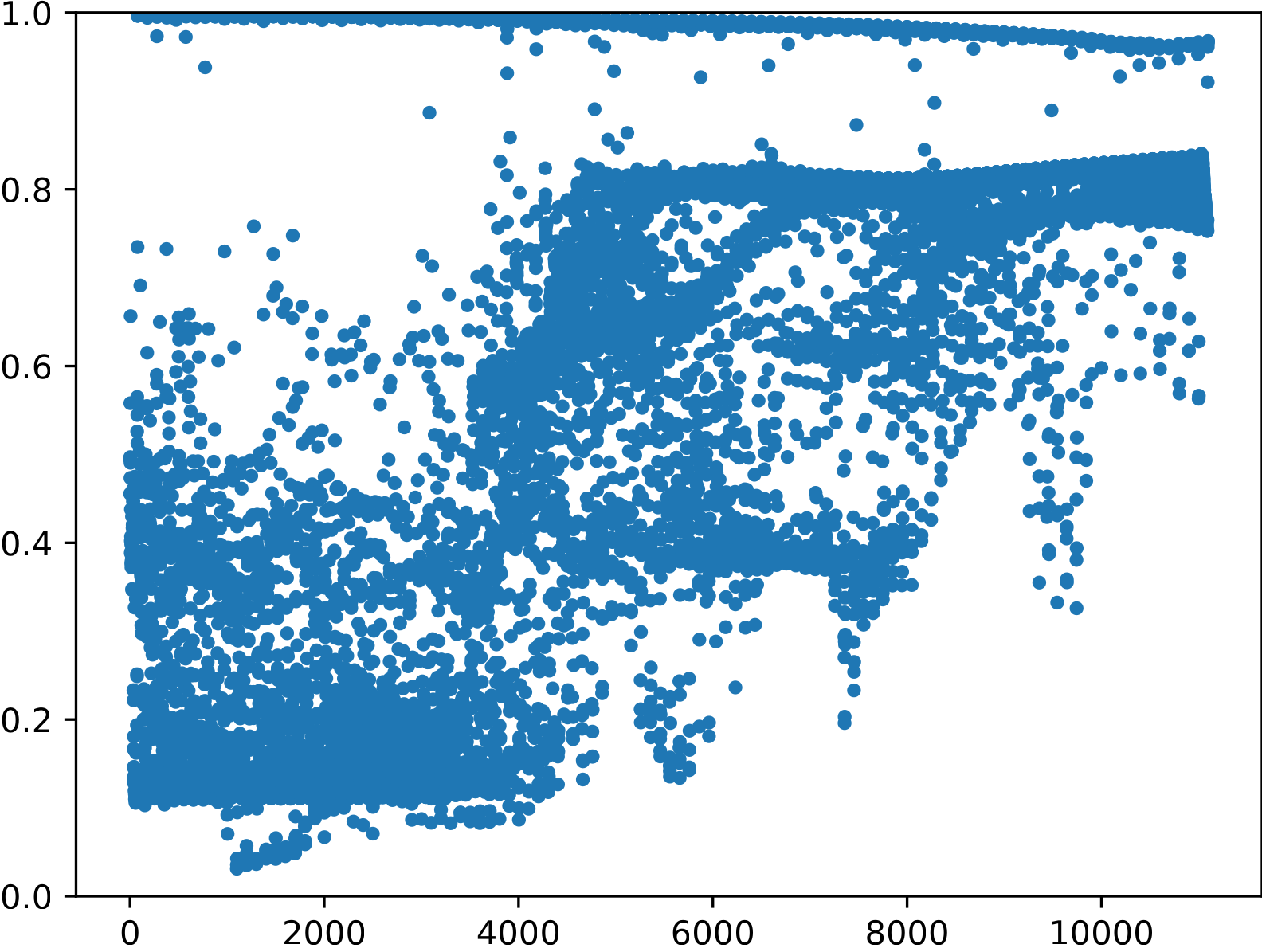} \\
		
		\ &
		\includegraphics[width=.3\textwidth]{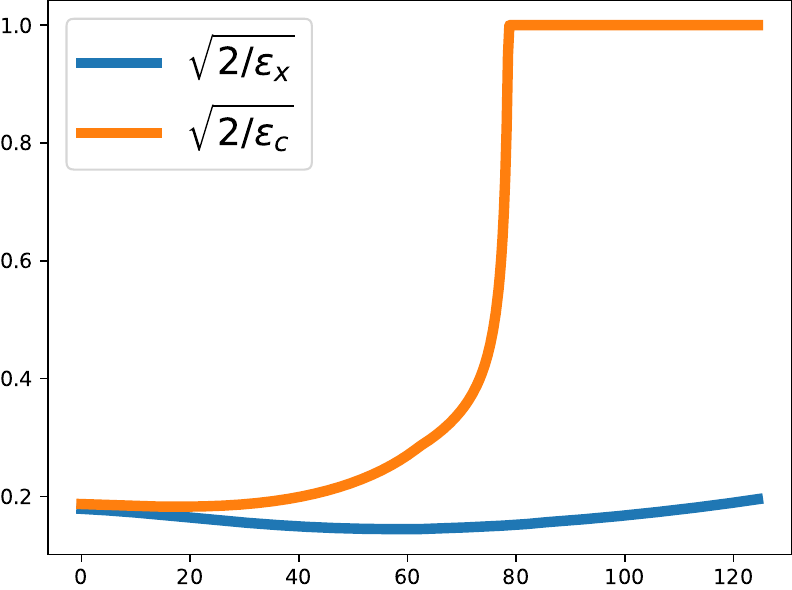} &
		\includegraphics[width=.3\textwidth]{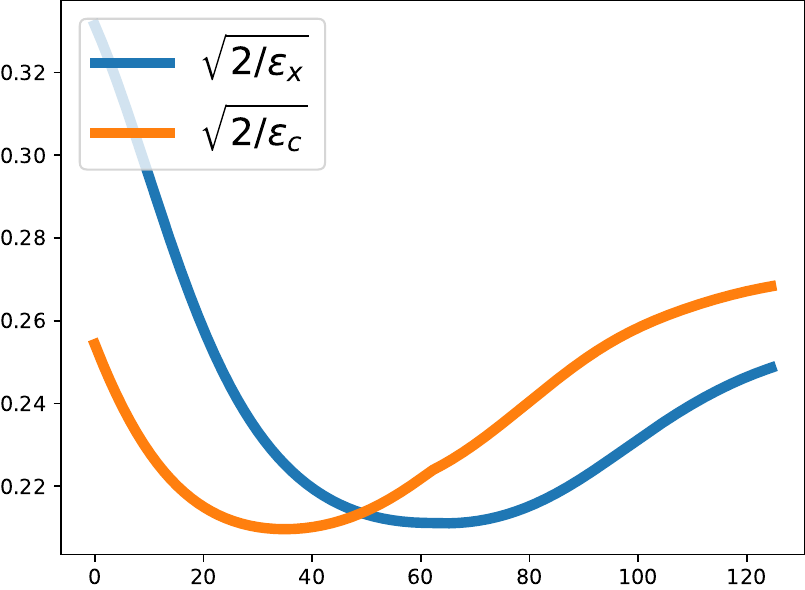}
	\end{tabular}
	\caption{Comparison of Algorithm \ref{alg:dal} and Algorithm \ref{alg:admm} on another test image portraying adult Namou. Output image on top, grey level for every pixel in the middle row, time evolution of the semi-axes of the ellipsoid \eqref{eq:ellipsoid} on the bottom.}
	\label{fig:namou2_admm}
\end{figure}

To conclude, we perform a final experiment on an image that presents mixed features. Looking at the original picture on the top-left of Figure \ref{fig:namou2_admm}, we notice both large areas with high contrast (the floor, the irregular wall on the right, the dark plants in the background) and a central subject with many details and shades of grey. 
The output of Algorithm \ref{alg:dal} is definitely unsatisfactory. Even if the larger areas are smoothed out, the diffusion effect due to global interactions in the colour component (see the bottom-centre plot in Figure \ref{fig:namou2_admm}) is dominant. As a consequence, not even the white wall is homogenised correctly and the railing in the background is almost fused with the trees. Overall, Algorithm \ref{alg:dal} acts more like a denoiser in this test. On the contrary, Algorithm \ref{alg:admm} preserves all the edges neatly, even enhancing the contrast in some areas. The floor and the wall on the side are evened out, and the dog is perfectly outlined. The difference between the two outputs can be seen in the colour distribution as well. While there are almost no clusters in the one given by Algorithm \ref{alg:dal}, we can clearly identify them in the middle-right plot. Indeed, the optimal controls found by the two algorithms are completely different, as can be seen on the bottom row of Figure \ref{fig:namou2_admm}.

\subsection{Extension to high-resolution images}

Once an optimal control is obtained with the ADMM algorithm, since by construction it is only related to the spatial configuration of the pixels, and not their number, it can seamlessly be applied to higher-resolution versions of the same image. This task only consists in one forward-in-time integration of the multi-agent system \eqref{eq:MAsystem}, which is much faster.

\begin{figure}[!h]
	\centering
	\begin{tabular}{cccc}
		\hline
		\textit{ADMM input} & \textit{ADMM output} & \textit{Hi-res image} & \textit{Quantised hi-res} \\ 
		\hline \\
		\includegraphics[width=2.7cm]{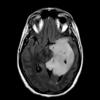} &
		\includegraphics[width=2.7cm]{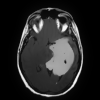} &
		\includegraphics[width=2.7cm]{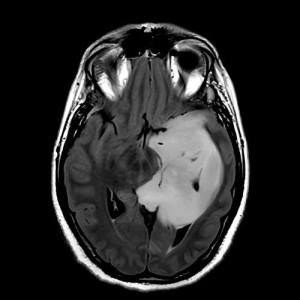} &
		\includegraphics[width=2.7cm]{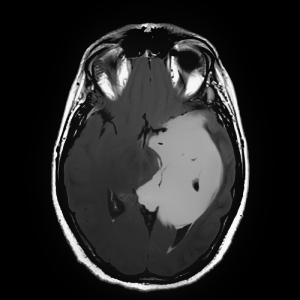} \\
		
		\includegraphics[width=2.7cm]{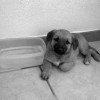} &
		\includegraphics[width=2.7cm]{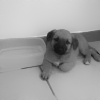} &
		\includegraphics[width=2.7cm]{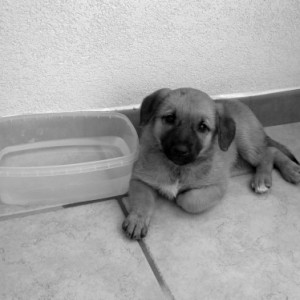} &
		\includegraphics[width=2.7cm]{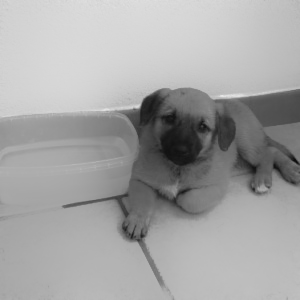} \\
		
		\includegraphics[width=2.7cm]{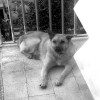} &
		\includegraphics[width=2.7cm,height=2.7cm]{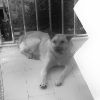} &
		\includegraphics[width=2.7cm]{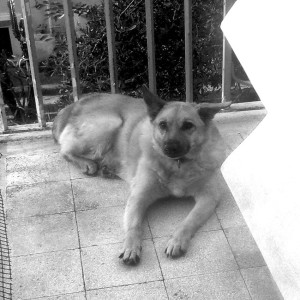} &
		\includegraphics[width=2.7cm]{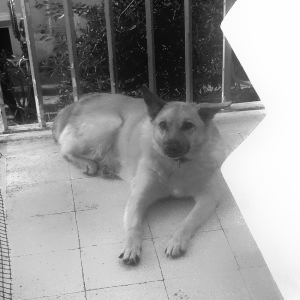}
	\end{tabular}
	\caption{Application of the controls obtained with Algorithm \ref{alg:admm} to high-resolution versions of the input images.}
	\label{fig : HD}
\end{figure}

As an example, we apply the controls obtained with Algorithm \ref{alg:admm} to a high-resolution version of the three test images from the previous sections. The results are shown in Figure \ref{fig : HD}. The output images are sharper and more detailed than the low-resolution versions, while still preserving the same segments. In particular, the tumour in the MRI scan is even more clearly outlined, and the dog in the second and third tests is perfectly defined, with a clear contrast with the background.

\section{Comparison with other methods}

In conclusion to this chapter, we present some comparisons between our model and other established tools for image denoising and segmentation. In particular, we consider the Chambolle--Pock (CP) algorithm for total variation minimisation~\cite{CP_2011}, the multi-agent model proposed by Herty--Pareschi--Visconti (HPV) in \cite{Herty2020}, and the non-local means (NLM) algorithm~\cite{buades2011non} by Buades--Coll--Morel. The first method is the most widely used for solving the ROF model, and it is based on a primal--dual formulation of the problem. The second one is a multi-agent system where the pixels are allowed to move during the evolution (but without any control variable) forming clusters not only in the colour feature, but also in space. The third one is a non-parametric method for image denoising, which relies on the idea of averaging similar patches in the image. We highlight that this last method is not designed for segmentation, but it is still interesting to compare it with our model, as it is based on a non-local interaction between pixels, which is one of the key features of our approach.

To perform a quantitative comparison, we use three metrics: intra-region variance, inter-region contrast, and a validity score that is the ratio between the previous two. Intra-region variance is computed following the approach introduced in \cite{levine1985dynamic}. First, using each output image from the aforementioned algorithms, the segments $R_1, \ldots, R_k$ are labeled with the Watershed algorithm~\cite{watershed}. Then, the labels are assigned to the original image and the following quantity is computed:
\begin{equation}\label{eq:intra_region_variance}
	\mathit{IRV} = \frac{1}{N} \sum_{j=1}^k |R_j| \, \sigma_j^2,
\end{equation}
where $|R_j|$ is the number of pixels in the segment $R_j$, and $\sigma_j^2$ is the variance of the pixel values in $R_j$ on the original image. The reason why this metric is evaluated on the original image, and not on the output of the algorithm, is that we want to assess how well the segments created by the algorithm are able to capture homogenous regions and objects in the original picture. In other words, we want to verify that the segments are not just uniform in colour, but that they are also meaningful: the lower the intra-region variance, the better the segmentation. However, $\mathit{IRV}$ alone is not sufficient to estimate the quality of the segmentation, as it can be minimised by creating a large number of small regions. In fact, a segmentation where each region contains a single pixel would achieve zero intra-region variance. To avoid this issue, we also compute the inter-region contrast as
\begin{equation}\label{eq:inter_region_contrast} 
	\mathit{IRC} = \min_{\substack{i,j=1,\ldots,k \\ i < j}} |m_i - m_j|,
\end{equation}
where $m_j$ is the mean colour of the segment $R_j$ on the processed image. This metric quantifies how well the segments are separated from each other in terms of colour. Over-segmentation is penalised by $\mathit{IRC}$, as it would yield a very small value, whereas a good segmentation should have a large inter-region contrast. Finally, to capture the overall quality of the segmentation, we define the validity score as
\begin{equation}\label{eq:validity_score}
	\mathit{VS} = \frac{\mathit{IRV}}{\mathit{IRC}},
\end{equation}
which should be as small as possible.

\begin{figure}
	\centering
	\begin{tabular}{c c c c c}
		\hline
		\textit{ADMM} & \textit{CP} & \textit{HPV} & \textit{NLM} \\
		\hline \\
		\includegraphics[width=2.4cm]{brain300_out.png} &
		\includegraphics[width=2.4cm]{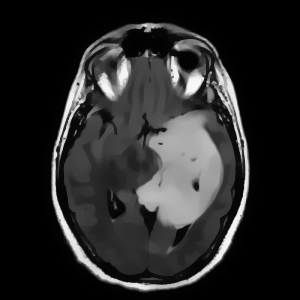} &
		\includegraphics[width=2.4cm]{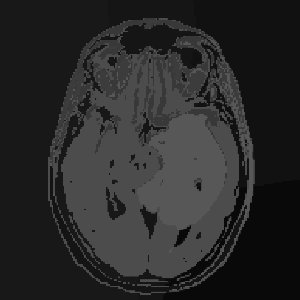} &
		\includegraphics[width=2.4cm]{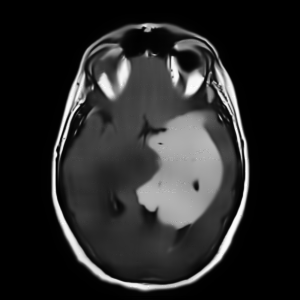} \\
		\includegraphics[width=2.4cm]{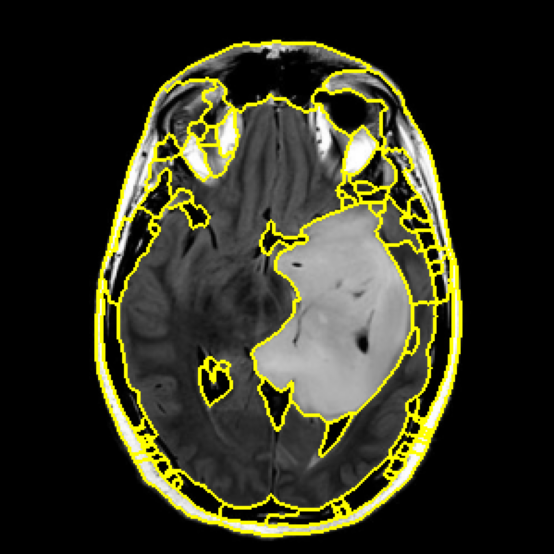} &
		\includegraphics[width=2.4cm]{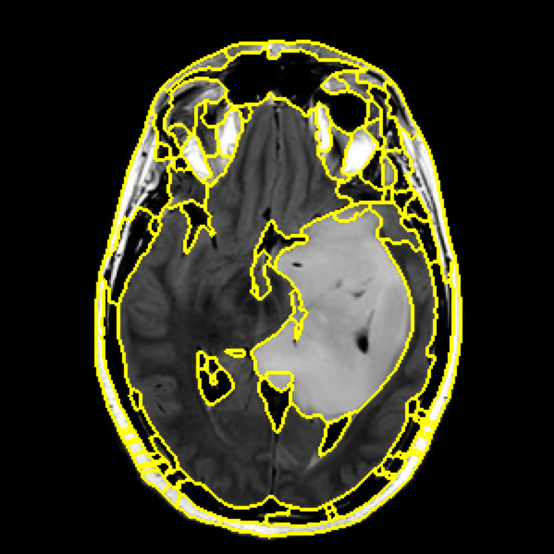} &
		\includegraphics[width=2.4cm]{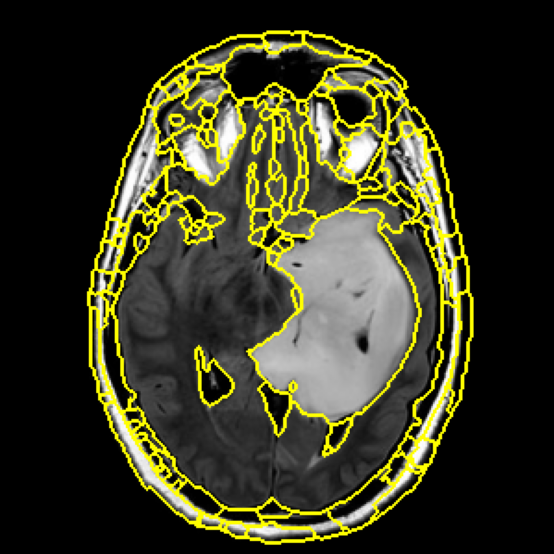} &
		\includegraphics[width=2.4cm]{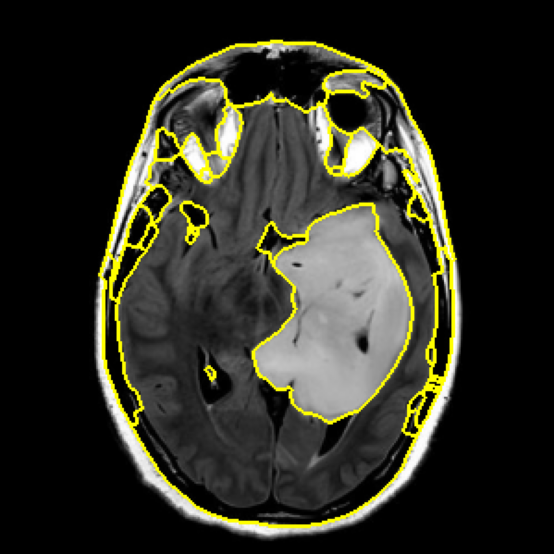} \\
		& & & & \\
		\includegraphics[width=2.4cm]{namou300_out.png} &
		\includegraphics[width=2.4cm]{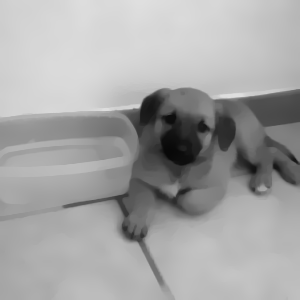} &
		\includegraphics[width=2.4cm]{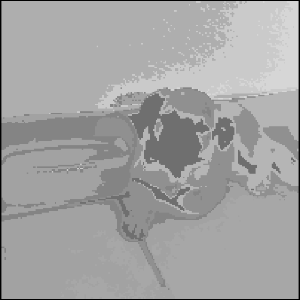} &
		\includegraphics[width=2.4cm]{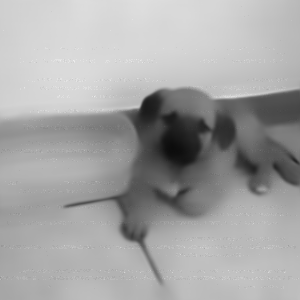} \\
		\includegraphics[width=2.4cm]{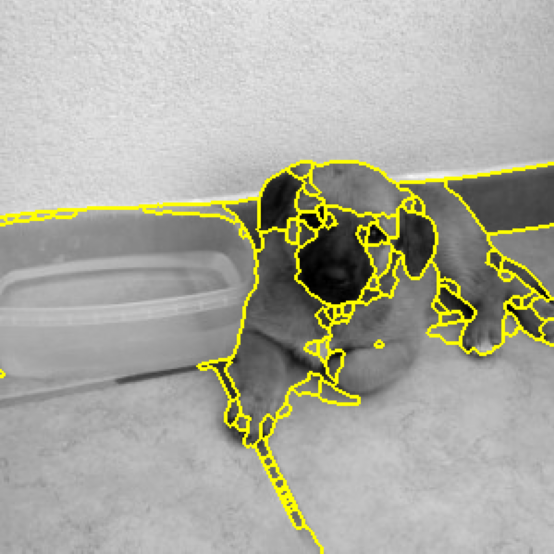} &
		\includegraphics[width=2.4cm]{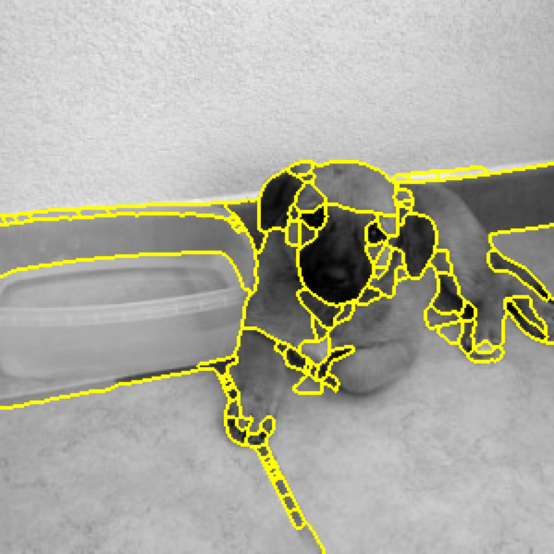} &
		\includegraphics[width=2.4cm]{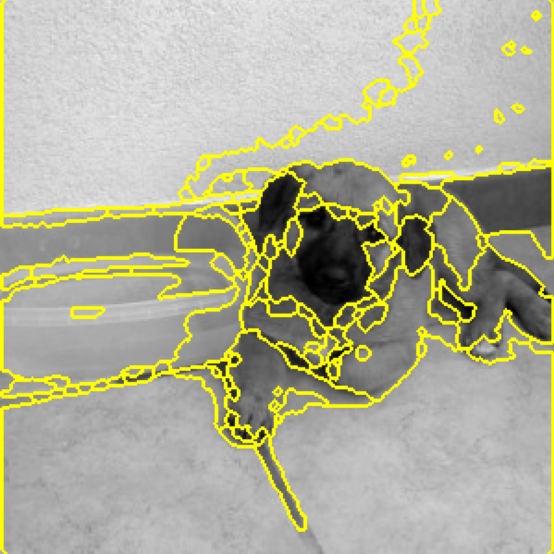} &
		\includegraphics[width=2.4cm]{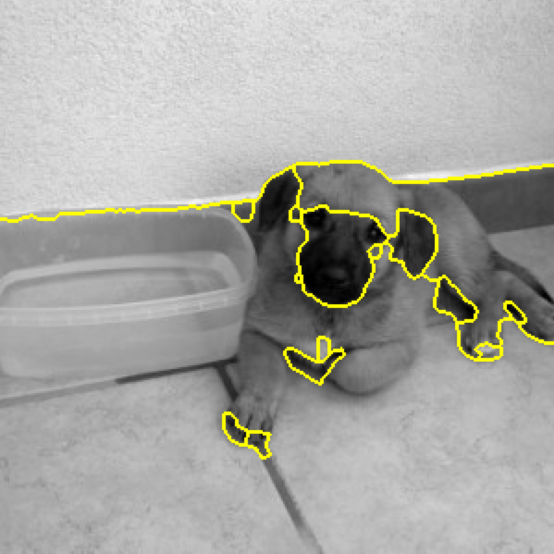} \\
		& & & & \\
		\includegraphics[width=2.4cm]{namou2300_out.png} &
		\includegraphics[width=2.4cm]{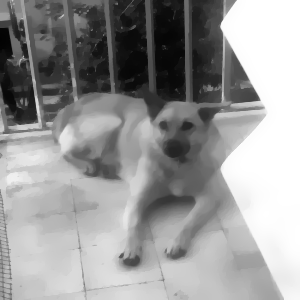} &
		\includegraphics[width=2.4cm]{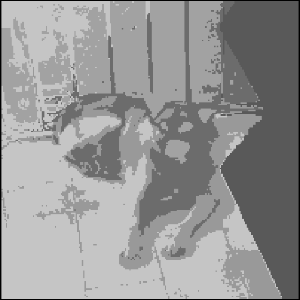} &
		\includegraphics[width=2.4cm]{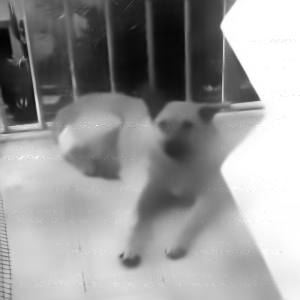} \\
		\includegraphics[width=2.4cm]{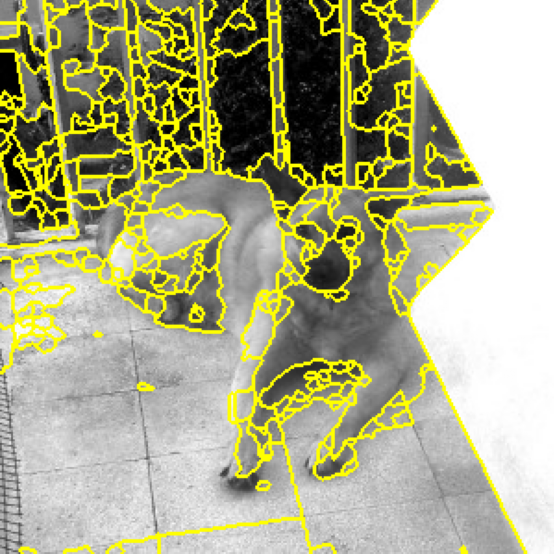} &
		\includegraphics[width=2.4cm]{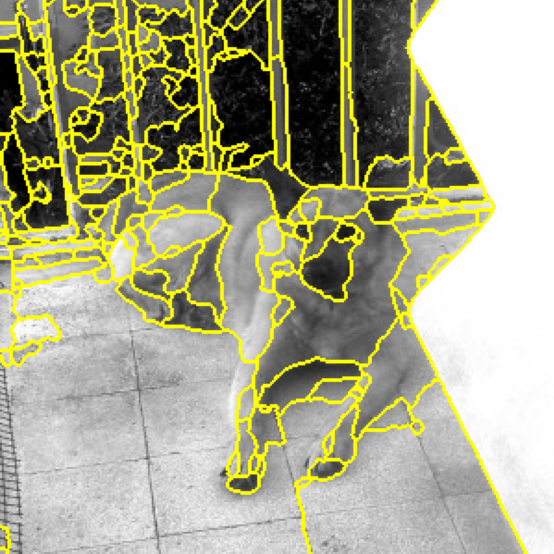} &
		\includegraphics[width=2.4cm]{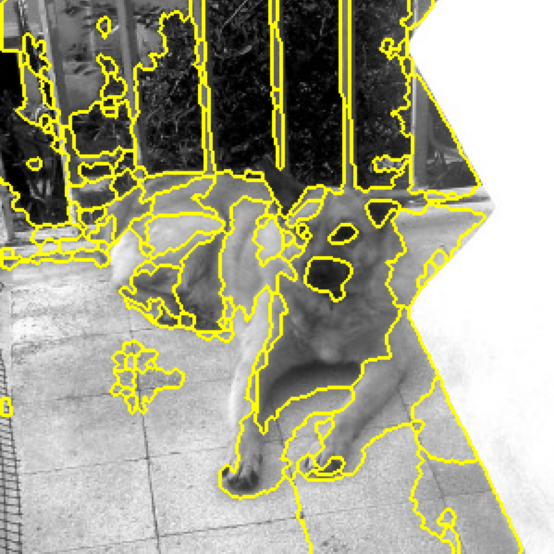} &
		\includegraphics[width=2.4cm]{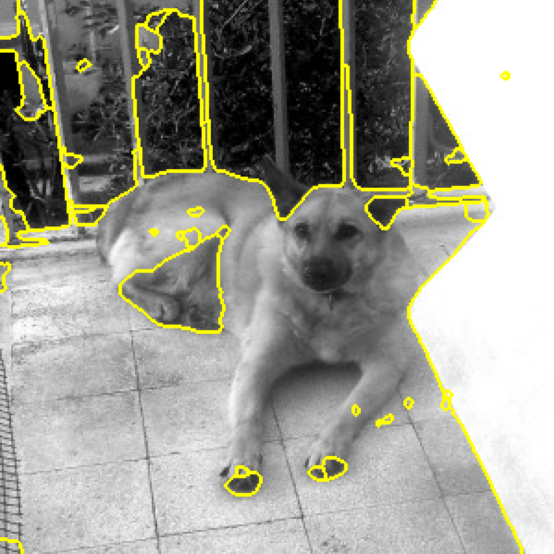}
	\end{tabular}
	\caption{Comparison of the outputs of Algorithm \ref{alg:admm}, the Chambolle--Pock algorithm, the multi-agent model by Herty--Pareschi--Visconti, and the non-local means algorithm on the test images from Figure \ref{fig : HD}. Under each picture, the regions identified by the Watershed algorithm are superimposed on the original image.}
	\label{fig:segmentation_comparison}
\end{figure}

We conduct the test on a $300\times 300$-pixel version of the three pictures used in the previous sections (see Figure \ref{fig : HD}). In order to get a fair comparison, we tune the parameters of the Chambolle--Pock algorithm so that the cost functional of the ROF model is the same as the one of Algorithm \ref{alg:admm} at convergence; in particular, we choose the fidelity parameter $\lambda = \alpha \, \Delta x$. For the multi-agent model by Herty et al., we set the bounded confidence values  $\epsilon_1=0.2$ and $\epsilon_2=0.1$. The threshold for the Watershed algorithm is set to $0.02$. The results are shown in Figure \ref{fig:segmentation_comparison} and the values of the three metrics are reported in Table \ref{tab:comparison}. Overall, the ADMM algorithm achieves very good validity scores, and it is competitive with the other methods. In particular, it is the best one in the first and third tests. Moreover, it is worth noting that the output of Algorithm \ref{alg:admm} is more visually appealing than the one of the other algorithms, as it preserves edges and details without over-segmenting the image. Finally, we remark that, despite apparently outperforming the other methods in Table \ref{tab:comparison}, the non-local means algorithm is not able to create meaningful segments in the images, as can be seen in Figure \ref{fig:segmentation_comparison}. As a matter of fact, the regions identified by the Watershed algorithm on its output are poorly related to the objects in the original images.

\begin{table}[h]
	\centering
	\begin{tabular}{c l r r r r}
		\toprule
		\textit{Image} & \multicolumn{1}{c}{\textit{Metric}} & 
		\multicolumn{1}{c}{\textit{ADMM}} & 
		\multicolumn{1}{c}{\textit{CP}} & 
		\multicolumn{1}{c}{\textit{HPV}} & 
		\multicolumn{1}{c}{\textit{NLM}} \\
		\midrule
		\multirow{3}{*}{
			\includegraphics[height=1.4cm]{brain300.png}
		}
		& $\mathit{IRV}$ $(\times 10^{-3})$ & $9.775$ & $6.793$ & $8.465$ & $18.872$ \\
		& $\mathit{IRC}$ $(\times 10^{-5})$ & $2.504$ & $0.750$ & $2.045$ & $3.357$ \\
		& $\mathit{VS}\phantom{0}$ $(\times 10^{2})$   & $3.904$ & $9.052$ & $4.139$ & $5.622$ \\
		
		\midrule
		
		\multirow{3}{*}{
			\includegraphics[height=1.4cm]{namou300.png}
		}
		& $\mathit{IRV}$ $(\times 10^{-3})$ & $5.920$ & $6.022$ & $2.547$ & $12.631$ \\
		& $\mathit{IRC}$ $(\times 10^{-5})$ & $3.961$ & $7.003$ & $1.548$ & $115.940$ \\
		& $\mathit{VS}\phantom{0}$ $(\times 10^{2})$   & $1.495$ & $0.860$ & $25.471$& $0.109$ \\
		
		\midrule
		
		\multirow{3}{*}{
			\includegraphics[height=1.4cm]{namou2300.png}
		}
		& $\mathit{IRV}$ $(\times 10^{-3})$ & $11.767$& $8.200$ & $19.981$& $27.652$ \\
		& $\mathit{IRC}$ $(\times 10^{-5})$ & $1.556$ & $0.100$ & $2.347$ & $4.312$ \\
		& $\mathit{VS}\phantom{0}$ $(\times 10^{2})$   & $7.561$ & $82.00$ 0& $29.520$& $6.414$ \\
		
		\bottomrule
	\end{tabular}
	\caption{Inter-region variance \eqref{eq:inter_region_contrast}, intra-region contrast \eqref{eq:inter_region_contrast} and validity score \eqref{eq:validity_score} for the outputs of Algorithm \ref{alg:admm}, the Chambolle--Pock algorithm, the Herty--Pareschi--Visconti multi-agent model and the non-local means algorithm on the test images from Figure \ref{fig : HD}. The leftmost column shows the original image for reference.}
	\label{tab:comparison}
\end{table}

\section{Conclusions and perspectives}
\label{sec:conclusions}

In this paper, we introduced a novel approach for image quantization and segmentation by framing them as an optimal control problem for a multi-agent system. The core idea is to control the colour dynamics of image pixels to form clusters that correspond to different segments of the image. This is achieved by defining a cost functional that balances colour variation and fidelity to the original image. Initially, a model based on a quadratic cost functional for the colour gradient was proposed and the optimisation problem was solved using a direct--adjoint looping (DAL) algorithm. Numerical experiments demonstrated that this approach is effective for images with high-contrast features, such as medical MRI scans. However, it was found to produce a blurring effect and was less suitable for natural images with smoother colour transitions, where it acted more like a denoising filter. To address these limitations, a more sophisticated model incorporating total variation was introduced. This model is based on the Rudin-Osher-Fatemi (ROF) functional, which is well-suited for preserving sharp edges. Due to the non-differentiability of the total variation term, the minimisation problem was regularised and reformulated as a saddle-point problem, which was solved numerically using the Alternating Direction Method of Multipliers (ADMM). This primal-dual approach effectively mitigates the diffusion and blurring effects observed with the first model. Comparative numerical tests on various images showcased the superior performance of the ADMM-based total variation model. It successfully segmented images with both high-contrast and complex, detailed regions, preserving sharp edges and creating more uniform colour clusters. The parallel implementation of the numerical schemes on a GPU using CUDA significantly improved computational performance, highlighting the potential of this methodology for handling high-dimensional data.

Future work could involve a comparison with other state-of-the-art segmentation tools, as well as explore the application of this framework to colour images. Further research could also investigate the mean-field limit of the model or introduce stochastic methods, in order to treat high-resolution images more efficiently.

\section*{Acknowledgments}
We declare that no dog was harmed for the execution of the numerical tests.

This work was supported by the Italian National Group of Scientific Computing (GNCS-INDAM). The authors also acknowledge the support of Sapienza University of Rome for the project ``Advanced Computational Methods for Real-World Applications: Data-Driven Models, Hyperbolic Equations and Optimal Control", CUP: B83C25000880005.

A.O. is supported by the European Union -- Next Generation EU,  Mission~4, Component 1, CUP 351: B83C22003530006.

S.C. is supported by the PNRR-MUR project ``Italian Research Center on High Performance Computing, Big Data and Quantum Computing''.

G.V. acknowledges the support of MUR (Ministry of University and Research) under the MUR-PRIN PNRR Project 2022 No.~P2022JC95T ``Data-driven discovery and control of multi-scale interacting artificial agent system''.

\bibliographystyle{abbrv}
\bibliography{references}

\end{document}